\documentclass[11pt,a4paper]{siamart220329}

\usepackage{braket}
\usepackage[utf8]{inputenc}
\usepackage[english]{babel}

\usepackage{tikz}
\usetikzlibrary{shapes.geometric}

\usepackage{graphicx}
\usepackage{fancyhdr}
\usepackage{dsfont}
\usepackage{stmaryrd}
\usepackage{bbm}

\usepackage{amsfonts}
\usepackage{amssymb}

\usepackage{amsmath}
\DeclareMathOperator*{\argmax}{arg\,max}

\usepackage{amsfonts}
\usepackage{amssymb}
\usepackage{mathtools}
\usepackage{enumitem}
\usepackage{multirow}

\usepackage{xcolor}

\usepackage{mathrsfs}

\usepackage{overpic} 

\usepackage{algorithm}
\usepackage[noend]{algpseudocode}

\usepackage{mathtools,booktabs}

\algrenewcommand\algorithmicrequire{\textbf{Input:}}
\algrenewcommand\algorithmicensure{\textbf{Output:}}

\usepackage{varwidth}

\usepackage[title]{appendix}

\newcommand{\rank}{{\rm rank}}

\usepackage{color}

\definecolor{rev1}{HTML}{cb270f}
\definecolor{rev2}{HTML}{1c8235}
\newcommand{\reva}[1]{{\color{black}{#1}}}
\newcommand{\revb}[1]{{\color{black}{#1}}}

\graphicspath{{images/}}

\usepackage[backend=bibtex]{biblatex}
\bibliography{references}

\title{Distributed memory parallel adaptive tensor-train cross approximation
\thanks{Submitted to the editors \today. \funding{T. Shi acknowledges support provided by the Director, Office of Science, Office of Advanced Scientific Computing Research , of the U.S. Department of Energy under Contract No. DE-AC02-05CH11231. D. Hayes and J.-M. Qiu acknowledge support provided by Department of Energy DE-SC0023164. J.-M. Qiu also acknowledges support provided by Air Force Office of Scientific Research FA9550-22-1-0390, FA9550-24-1-0254 and National Science Foundation NSF-DMS-2111253.}}}
\author{Tianyi Shi\thanks{Scalable Solvers Group, Lawrence Berkeley National Laboratory, Berkeley, CA 94720. (\email{tianyishi@lbl.gov})}
\and Daniel Hayes\thanks{Department of Mathematical Sciences, University of Delaware, Newark, DE 19716. (\email{dphayes@udel.edu})}
\and Jingmei Qiu\thanks{Department of Mathematical Sciences, University of Delaware, Newark, DE 19716. (\email{jingqiu@udel.edu})}}

\headers{Distributed memory Tensor-train cross}{T. Shi, D. Hayes, and J. Qiu}

\begin{document}
\newcommand{\R}[0]{\mathbb{R}}
\newcommand{\C}[0]{\mathbb{C}}
\maketitle

\begin{abstract}
The tensor-train (TT) format is a data-sparse tensor representation commonly used in high dimensional function approximations arising from computational and data sciences. Various sequential and parallel TT decomposition algorithms have been proposed for different tensor inputs and assumptions. In this paper, we propose subtensor parallel adaptive TT cross, which partitions a tensor onto distributed memory machines with multidimensional process grids, and constructs an TT approximation iteratively with tensor elements. We derive two iterative formulations for pivot selection and TT core construction under the distributed memory setting, conduct communication and scaling analysis of the algorithm, and illustrate its performance with multiple test experiments. These include up to 6D Hilbert tensors and tensors constructed from Maxwellian distribution functions that arise in kinetic theory. Our results demonstrate significant accuracy with greatly reduced storage requirements via the TT cross approximation. Furthermore, we demonstrate good to optimal strong and weak scaling performance for the proposed parallel algorithm.
\end{abstract}



\section{Introduction}\label{sec:introduction}
The success and development of computing machines in the past few decades have allowed researchers to deal with high-dimensional data easily. Typically, these data sets are stored as multidimensional arrays called tensors~\cite{kolda2009tensor}, and a general tensor $\mathcal{X} \in \C^{n_1\times\cdots \times n_d}$ requires a storage cost of $\prod_{j=1}^d n_j$ degrees of freedom. This scales exponentially with the dimension $d$, and is often referred to as ``the curse of dimensionality". Therefore, data-sparse tensor formats such as canonical polyadic (CP)~\cite{hitchcock1927expression}, Tucker~\cite{de2000multilinear}, tensor-train (TT)~\cite{oseledets2011tensor}, and tensor networks~\cite{evenbly2011tensor} with more complex geometries have been proposed. In particular, the TT format, also known as the matrix product state (MPS) in tensor networks and quantum physics, has a memory footprint that scales linearly with respect to the mode sizes $n_j$ and the dimension $d$. The TT format is widely used in applications such as molecular simulations~\cite{savostyanov2014exact}, high-order correlation functions~\cite{kressner2015low}, partial differential equations~\cite{guo2023local}, constrained optimization~\cite{dolgov2017low,benner2020low}, and machine learning~\cite{vandereycken2022ttml,novikov2020tensor}. Furthermore, the TT format can be incorporated with extra conditions to form special tensor representations that can capture latent data structures. For example, the quantized TT~\cite{dolgov2012fast} format is a combination of the TT format and hierarchical structures, and the tensor chain format ~\cite{espig2012note} is a result of alterations on MPS.

In practice, instead of finding an exact TT representation of a tensor $\mathcal{X}$, one aims to construct an approximation $\tilde{\mathcal{X}}$ with a low rank format. One group of TT decomposition algorithms targets approximations with a controllable error estimate,
\begin{equation}
\| \mathcal{X} - \tilde{\mathcal{X}} \|_F \leq \epsilon \| \mathcal{X} \|_F, \qquad \|\mathcal{X}\|_F^2 = \sum_{i_1=1}^{n_1} \cdots \sum_{i_d = 1}^{n_d} |\mathcal{X}_{i_1,\ldots,i_d}|^2,
\label{eq:FrobeniusNorm}
\end{equation}
where $0\leq\epsilon<1$ is an accuracy tolerance~\cite{grasedyck2013literature,hackbusch2012tensor} and the bounds hold with stability and high probability. Such algorithms include TT singular value decomposition (TTSVD)~\cite{oseledets2011tensor}, TT sketching~\cite{che2019randomized}, and QR or SVD based TT cross~\cite{oseledets2010tt,dektor2024collocation}.
Another category of TT approximation algorithms is to build $\tilde{\mathcal{X}}$ iteratively with greedy and heuristic approaches, such as TT alternating least squares~\cite{holtz2012alternating} and adaptive TT cross~\cite{oseledets2010tt,dolgov2020parallel}. Although we do not have theoretical guarantees for convergence or convergence rates, these methods can have good performance in certain scenarios. Particularly, adaptive TT cross is a data-based algorithm, with small complexity cost and especially suitable for sparse, giant datasets.

In order to exploit modern computing architectures, researchers have proposed various parallel methods, including both shared memory and distributed memory parallelism, for tensor decomposition in CP~\cite{li2017model,smith2015splatt}, Tucker and hierarchical Tucker~\cite{austin2016parallel,grasedyck2019parallel,ballard2020tuckermpi,kaya2016high}, and TT~\cite{shi2023parallel,grigori2020parallel,chen2017parallelized,wang2020adtt,dolgov2020parallel} format. In addition, tensor operations in TT format, including addition and multiplication~\cite{daas2022parallel}, contraction~\cite{solomonik2014massively}, and recompression~\cite{al2023randomized} can be executed in parallel as well. In particular, distributed memory parallelism is an operational vessel for data parallel algorithms, so it is not compatible with tensor algorithms that require all elements at once, such as SVD. The major challenge we address in this paper is to construct an approximation $\tilde{\mathcal{X}}$ in a TT format from large $\mathcal{X}$ with distributed memory parallelism. This allows one to partition $\mathcal{X}$ into smaller blocks so that each process handles a significantly smaller chunk of data. Furthermore, with a successful distributed memory design, all processes can execute individual shared memory parallel algorithms and minimize communications with each other, leading to efficient computational and storage consumption. 

In this paper, we propose a new parallel adaptive TT cross algorithm based on a distributed memory framework built upon subtensors. \reva{This is viewed as a direct parallelization of the construction in~\cite[Section 3]{oseledets2010tt}.} A subtensor is a multilinear generalization of a submatrix, and has been used in the matricized-tensor times Khatri-Rao product (MTTKRP)~\cite{ballard2018, ballard2020}, hierarchical subtensor decomposition~\cite{ehrlacher2021}, parallel Tucker decomposition~\cite{ballard2020tuckermpi}, and parallel TT decomposition~\cite{shi2023parallel}. The subtensors can carry out key kernels in adaptive TT cross independently. In the end, results on subtensors are gathered to form the outcome of the entire tensor. For the remainder of this paper, we call this subtensor parallelism, as opposed to dimension parallelism. In fact, one can understand subtensor parallelism as a special type of data parallelism, as the key process is to distribute tensor elements to computing architectures in a regular pattern. Comparatively, dimension parallel algorithms partition computations with respect to the dimensionality of the tensor, and thus the number of processes used actively is limited even in the distributed memory setting. The processes can also encounter severe load imbalance as computations for each dimension may vary significantly. In subtensor parallelism, we construct a multidimensional process grid for subtensor partitioning, which enables us to derive explicit bounds on the bandwidth and communication costs. As a bonus, we can run dimension parallel algorithms in the shared memory setting on each process.

In many applications such as numerical integration in quantum mechanics~\cite{meyer1990multi} and inverse problems with uncertainty~\cite{stuart2010inverse}, and data analysis in statistical science~\cite{mccullagh2018tensor} and machine learning~\cite{rabanser2017introduction}, tensors are often formed without an exact formula and can be extremely sparse. In these cases, researchers develop adaptive TT cross~\cite{oseledets2010tt} and dimension parallel adaptive TT cross~\cite{dolgov2020parallel} for data-centric TT approximation. Our main contribution in this paper is to develop adaptive matrix and tensor cross approximations within the distributed memory framework using submatrices and subtensors. In particular, we derive two novel communication-efficient iterative procedures to construct matrix cross approximations and show they can recover accurate results. These two iterative procedures are combined with submatrix cross, and are used for TT pivot selection and TT core construction respectively. Especially, we show that local pivots selected on the subtensors still hold the nestedness property to maintain tensor interpolation for all steps and all dimensionalities. Furthermore, we can apply dimension parallel TT core construction on each subtensor in a shared memory setting, and achieve a comparably good approximation for the entire tensor.

We implement our parallel algorithms with both distributed and shared memory parallelism in Python. Particularly, we use MPI4Py for distributed memory setup, which is a Python framework of the message passing interface (MPI). In addition, we use numpy for linear algebra operations to optimize our codes, which is a Python wrapper for well-established linear algebra packages such as BLAS and LAPACK. 
The remaining of the manuscript is organized as follows.~\Cref{sec:background} reviews some necessary tensor notations and the TT format with existing serial and parallel algorithms. In~\cref{sec:subTTcross}, we introduce the new subtensor parallel TT cross algorithm. Then, we provide scalability and complexity analysis in~\cref{sec:subComm}. Finally, we demonstrate the performance on up to 6D
datasets in~\cref{sec:NumericalExamples}.

\section{Tensor notations, TT format, and TT cross} \label{sec:background}
In this section, we review some tensor notations, the TT format for low rank tensor approximations, and the adaptive TT cross approximation.

\subsection{Tensor notation} \label{sec:notation}
We use lower case letters for vectors, capital letters for matrices, and calligraphic capital letters for tensors. \revb{Just as tensors are higher order analogues of matrices, subtensors are also higher order analogues of submatrices, storing parts of an entire tensor.} For notational simplicity, we use MATLAB-style notation to start index counting from 1, and the symbol ``:" in indexing. This includes using $a\!:\!b$ to represent the inclusive index set $\{a,a+1,\ldots,b\}$, and a single ``:" to represent all the indices in that dimension from start to end. For example, if $\mathcal{Y} \in \R^{4 \times 8 \times 10}$, then $\mathcal{Y}(\ :\ ,4\!:\!6,\ :\ )$ or $\mathcal{Y}_{:,4:6,:}$ denotes a subtensor of $\mathcal{Y}$ with size $4 \times 3 \times 10$. In addition, we use index sets for submatrix and subtensor selection. For example, $A(:,\mathcal{J})$ is a submatrix of $A$ with columns $A_{:,j}$ for all $j \in \mathcal{J}$.

Furthermore, we use the MATLAB command ``reshape" to form a new structure according to the multi-index via reorganizing elements without changing the element ordering. For example, if $\mathcal{Y} \in \C^{n_1 \times n_2 \times n_3}$, then $Z = {\rm reshape}(\mathcal{Y},n_1n_2,n_3)$ returns a matrix of size $n_1n_2 \times n_3$ formed by stacking entries, and similarly, $\mathcal{Y} = {\rm reshape}(Z,n_1,n_2,n_3)$. The command ``reshape" is essential when flattening a tensor into matrices, which we refer to as the unfoldings of a tensor. Tensor unfoldings are fundamental to the TT format, especially in developing decomposition algorithms and bounding TT ranks. For a tensor $\mathcal{X}\in\C^{n_1\times\cdots\times n_d}$, we denote the $k$th unfolding as
\[X_k={\rm reshape}\left(\mathcal{X},\prod_{s=1}^k n_s,\prod_{s=k+1}^d n_s\right), \quad 1 \le k \le d-1.\]

\begingroup
\setlength{\tabcolsep}{8pt} 
\renewcommand{\arraystretch}{1.5} 
\begin{table}[]
\reva{
\caption{Notation conventions for tensors and index sets.}
\label{tbl:notations}
\centering
\begin{tabular}{cc}
\hline
\hline
Concepts & Notations \\
\hline
Vector & Lower case letter $u,v,s$ \\
Matrix & Capital letter $A, C, U, R$ \\
Tensor & Calligraphic capital letter $\mathcal{X}, \mathcal{G}, \mathcal{T}$ \\
Single index & Lower case letter $i, j, k, \ell$ \\
Index set & Calligraphic capital letter $\mathcal{I}, \mathcal{J}, \mathcal{K}, \mathcal{L}$ \\
Index set with all possible indices & Blackboard bold capital letter $\mathbb{I}, \mathbb{J}$ \\
$k$-th Tensor unfolding & Capital letter with subscript $X_k$ \\
\hline
\hline
\end{tabular}}
\end{table}
\endgroup

\subsection{Tensor-train format} \label{sec:TT}
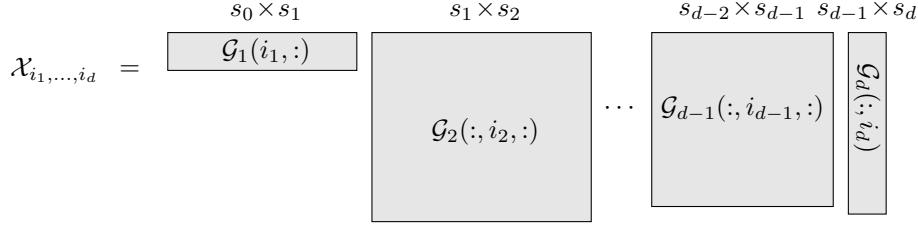
\begin{figure}
\centering
\begin{tikzpicture}
\filldraw[black] (0,-0.5) node {$\mathcal{X}_{i_1,\ldots,i_d}$};
\filldraw[black] (1,-0.5) node {$=$};
\filldraw[color=black,fill=gray!20] (1.5,0) rectangle (4,-.5);
\filldraw[black] (2.8,-0.25) node {$\mathcal{G}_1(i_1,:)$};
\filldraw[black] (2.8,0.3) node {$s_0 \! \times \! s_1$};
\filldraw[color=black,fill=gray!20] (4.2,0) rectangle (7.1,-2.5);
\filldraw[black] (5.7,-1.3) node {$\mathcal{G}_2(:,i_2,:)$};
\filldraw[black] (5.7,0.3) node {$s_1 \! \times \! s_2$};
\filldraw[black] (7.5,-1) node {$\cdots$};
\filldraw[color=black,fill=gray!20] (7.9,0) rectangle (10.3,-2.3);
 \filldraw[black] (9.1,-1) node {$\mathcal{G}_{d-1}(:,i_{d-1},:)$};
\filldraw[black] (9.1,0.3) node {$s_{d-2} \! \times \! s_{d-1}$};
\filldraw[color=black,fill=gray!20] (10.5,0) rectangle (11,-2.4);
\filldraw[black] (10.75,-1) node {\rotatebox{270}{$\mathcal{G}_{d}(:,i_{d})$}};
\filldraw[black] (10.75,0.3) node {$s_{d-1} \! \times \! s_d$};
\end{tikzpicture}
\caption{The TT format with TT core size $\pmb{s} = (s_0,\ldots,s_d)$. Each entry of a tensor is represented by the product of $d$ matrices, where the $k$th matrix in the ``train" is selected based on the value of $i_k$.}
\label{fig:TT}
\end{figure}

The TT format of a tensor $\mathcal{X}\in\C^{n_1\times \cdots \times n_d}$ comprises of $d$ TT cores, $\mathcal{G}_k \in \C^{s_{k-1} \times n_k \times s_k}$ for $1 \le k \le d$, and takes the representation
\[
\mathcal{X}_{i_1,\ldots,i_d} = \mathcal{G}_1(:,i_1,:)\mathcal{G}_2(:,i_2,:) \cdots \mathcal{G}_d(:,i_d,:), \qquad 1\leq i_k \leq n_k.
\]
In other words, each element in $\mathcal{X}$ can be computed as the product of a sequence of matrices. The vector $\pmb{s} = (s_0,\ldots,s_d)$ is referred to as the size of the TT cores, and in order for the product to be compatible, we require $s_0 = s_d = 1$. This TT representation thus has a storage cost of $\sum_{k=1}^d s_{k-1}s_k n_k$, which is linear with respect to both $d$ and $n_k$. In addition, we call a vector $\pmb{r}= (r_0,\ldots,r_d)$ the TT rank if $\pmb{r}$ contains entry-by-entry smallest possible values of the TT core size~\cite{oseledets2011tensor,shi2021compressibility}. In practice, the exact TT rank is hard to recover, so we either aim to obtain quasi-optimal TT core size from a given threshold, or build an accurate tensor approximation with small pre-selected TT core size. It is shown in~\cite{oseledets2011tensor} that ranks of tensor unfoldings bound the TT rank from above, so we hope to use $\pmb{s}$ that satisfies
\begin{equation} \label{eq:TT_trivial}
r_k \le s_k \le {\rm rank}(X_k), \quad 1 \le k \le d-1,
\end{equation}
where $\rank(X_k)$ is the rank of the $k$th unfolding of $\mathcal{X}$. In this way, if the ranks of all $X_k$ are small, the TT format is a memory-efficient representation.~\cref{fig:TT} provides an illustration of a tensor $\mathcal{X}$ in the TT format with TT core size $\pmb{s}$, and one may visualize slices of the TT cores as ``trains".

\subsection{Pivot selection and TT cross}
A major group of algorithms for TT decomposition is based on singular value decomposition (SVD) and randomized SVD, with deterministic or probabilistic error analysis. In practice, however, these SVD-based algorithms can be extremely expensive and unnecessary. For example, tensors originated from real world datasets are often large and sparse, and we often only want an approximation with an accuracy of a few significant digits. This leads to heuristic data-based methods, such as adaptive cross approximation (ACA). In summary, ACA mainly finds a ``skeleton" of original data for approximation. For example, \revb{using ACA, one can build the CUR factorization of a matrix, which is in the form of}
\begin{equation} \label{eq:matrix_CUR} 
A \approx CUR = A(:,\mathcal{J})A(\mathcal{I},\mathcal{J})^{-1}A(\mathcal{I},:), 
\end{equation}
where $\mathcal{I}$ and $\mathcal{J}$ are two index sets. The key building block in ACA is new ``pivot" selection to iteratively enrich these index sets. \revb{The selected indices to represent the data in ACA are referred to as ``pivots" due to their connection with the LU factorization. In fact, as shown later in the paragraph, the greedy approach to obtain ``pivots" is the same as performing multiple rounds of partial pivoting in Gaussian elimination.} 
Matrix ACA can be extended naturally to adaptive TT cross~\cite{oseledets2010tt}. For conceptual simplicity, all cross approximation algorithms we discuss in this article thereafter are with respect to adaptive cross approximations. There are various metrics for pivot selection, but we focus on the greedy approach outlined in~\cite[Algorithm 2]{dolgov2020parallel} \reva{(we include as~\cref{alg:MatCross} for self-containment)}, which aims to find an entry with the largest difference between the current approximation and the actual value. This method finds a quasi-optimal pivot \reva{through heuristics, and thus there is not a theoretical guarantee for accuracy,} but is computationally much cheaper than other routines such as max volume selection. Suppose we use $\mathcal{I}_{\le k}$ and $\mathcal{J}_{>k}$ to denote the selected sets of pivots for \revb{the row and column indices of $X_k$,} the $k$th unfolding of a tensor $\mathcal{X}$ with $1 \le k \le d-1$, and $\mathbb{I}_\ell$ to denote the set of all indices for dimension $\ell$ with $1 \le \ell \le d$, then one step of finding new pivots is described in~\cref{alg:TTcross}. Notice that $\mathcal{I}_{\le 0}$ and $\mathcal{J}_{>d}$ here are empty sets\revb{, entries of $\mathcal{I}_{\le k}$ fall between 1 and $\prod_{j=1}^kn_j$, and entries of $\mathcal{J}_{>k}$ fall between 1 and $\prod_{j=k+1}^dn_j$. For a pictorial illustration of the keyword ``superblock" used in~\cref{alg:TTcross}, we refer the readers to~\cite{dolgov2020parallel}.} Additionally, pivots selected this way can be shown to satisfy the nestedness property, which ensures that pivots found in one tensor unfolding can be carried over to subsequent unfoldings. As a result,~\cref{alg:TTcross} is automatically a dimension parallel algorithm.

\reva{
\begin{algorithm}
\caption{One step of finding new pivots in matrix cross.}
\begin{algorithmic}[1]
\label{alg:MatCross}
\reva{\Require{A matrix $X \in 
\C^{n_1 \times n_2}$, index sets $(\mathcal{I},\mathcal{J})$ containing current pivots, and current approximation $\tilde{X}$}
\Ensure{New pivots$(i^*, j^*)$.}
\State Pick a random set of samples $\mathcal{L}$ and select $(i^*,j^*)\leftarrow \argmax_{(i,j)\in\mathcal{L}} |X(i,j)-\tilde{X}(i,j)|$.
\While {$|X(i^*,j^*)-\tilde{X}(i^*,j^*)| < |X(i,j)-\tilde{X}(i,j)|$ for all $(i,j)$ such that $i=i^*$ or $j=j^*$}
\State Select $(i^*,j^*)\leftarrow \argmax_{i \in \mathbb{I}} |X(i,j^*)-\tilde{X}(i,j^*)|$.
\State Select $(i^*,j^*)\leftarrow \argmax_{j \in \mathbb{J}} |X(i^*,j)-\tilde{X}(i^*,j)|$.
\EndWhile}
\end{algorithmic}
\end{algorithm}}

\begin{algorithm}
\caption{One step of finding new pivots in TT cross.}
\begin{algorithmic}[1]
\label{alg:TTcross}
\Require {A tensor $\mathcal{X} \in \C^{n_1 \times \dots \times n_d}$, index sets $(\mathcal{I}_{\le k}, \mathcal{J}_{>k})$ containing current pivots \revb{for $1 \le k \le d-1$}.}
\Ensure {New pivots $(i^*_{\le k}, j^*_{>k})$ \revb{for $1 \le k \le d-1$.}}
\For {$1 \le k \le d-1$}
\State Use~\reva{\cref{alg:MatCross}} to find $(i^*_{\le k}, j^*_{>k})$ on a superblock $X_k(\mathcal{I}_{\le k-1}\times\mathbb{I}_k,\mathbb{I}_{k+1}\times \mathcal{J}_{> k+1})$, which is seen as a submatrix of the $k$th unfolding of $\mathcal{X}$.
\EndFor
\end{algorithmic}
\end{algorithm}

Finally, once all pivots are collected in the index sets $(\mathcal{I}_{\le k}, \mathcal{J}_{>k})$ for $1 \le k \le d-1$, then the tensor approximation $\tilde{\mathcal{X}}$ can be built element-wise as
\begin{align}
    \tilde{\mathcal{X}}(i_1,\dots,i_d) &= \mathcal{X}(i_1,\mathcal{J}_{>1})\left[\mathcal{X}(\mathcal{I}_{\le 1},\mathcal{J}_{>1})\right]^{-1}\mathcal{X}(\mathcal{I}_{\le 1},i_2,\mathcal{J}_{>2}) \nonumber \\
    &\times \left[\mathcal{X}(\mathcal{I}_{\le 2},\mathcal{J}_{>2})\right]^{-1} \cdots \mathcal{X}(\mathcal{I}_{\le d-1},i_d), \label{eq:ttcross_build_elem}
\end{align}
where the element access operator $()$ are overloaded for notational simplicity\revb{, and all inverses are matrix inverses. The use of overloading allows for consistent notation to avoid mixing terms to select 3-way tensors with $\mathcal{X}(\cdot,\cdot,\cdot)$ and matrices with $X_k(\cdot,\cdot)$.} \reva{This can be seen as an extension of~\cref{eq:matrix_CUR} in the TT format, with multiple subtensor selections and cross section inverses.} In this way, \reva{we can group consecutive terms in the previous expression to construct the TT cores as}
\begin{align}
    \mathcal{G}_1 &= \mathcal{X}(:,\mathcal{J}_{>1})\left[\mathcal{X}(\mathcal{I}_{\le 1},\mathcal{J}_{>1})\right]^{-1}, \ \mathcal{G}_d = \mathcal{X}(\mathcal{I}_{\le d-1},:), \nonumber \\
    \mathcal{G}_{k} = &\mathcal{X}(\mathcal{I}_{\le k-1},:,\mathcal{J}_{>k})\left[\mathcal{X}(\mathcal{I}_{\le k},\mathcal{J}_{>k})\right]^{-1}, \ 2 \le k \le d-1.\label{eq:ttcross_build}
\end{align}

\section{Subtensor parallelism for TT cross}
\label{sec:subTTcross}
In this section, we develop subtensor parallel TT cross suitable for the distributed memory framework. \revb{In particular, the distributed computing processes handle non-overlapping subtensors, but work together to obtain the TT format of the entire tensor.} Throughout this section, we assume each mode size $n_j$ is partitioned evenly into $C_j$ pieces for $1 \le j \le d$. In this way, there are $C = \prod_{j=1}^d C_j$ subtensors in total with roughly the same size. For the moment, we suppose one process handles one subtensor at a time. This assumption can easily be lifted so that the subtensor grid and the process grid are totally different. In the special case that $d=1$ or $d=2$, these grids are often referred to as 1D or 2D grid respectively in literature~\cite{blackford1997scalapack}. In order to clearly describe our distributed memory algorithm, we refer to some MPI terminology and functions for communication patterns. These include
\begin{itemize}[leftmargin=*,noitemsep]
\item \textit{Root}: one process in a communication group to initialize collective operations.
\item \textit{Send}: The action of sending some information from one MPI rank to another.
\item \textit{Receive}: The action of receiving the information sent by \textit{Send}.
\item \textit{Gather}: The action of collecting some information from all processes to the root.
\item \textit{Allgather}: Same as \textit{Gather} but the information is stored on each process.
\end{itemize}

\subsection{Submatrix parallel matrix cross approximation}
\label{sec:submatrixCross}
We begin by introducing submatrix parallel matrix cross, and later use it to develop subtensor parallel TT cross. To start, we develop a new iterative formulation to construct ACA.

\subsubsection{A \revb{derivation of an} iterative construction of matrix cross approximation}
\label{sec:iterMat1}
\revb{In the following section, we go through the derivation of an iterative construction of matrix cross approximations. The final formula given in~\cref{eq:aca_approx_E} is originally proposed in~\cite[Section 2]{bebendorf2000approximation}, but we show a full proof for completeness. In addition,~\cref{eq:aca_approx_T} in~\cref{sec:iterMat2} is a corollary of the derivation, and to the best of our knowledge, it is a new formula for core construction in the tensor setting.}

Assuming that at step $z$ in ACA, we already have row and column indices $\mathcal{I}, \mathcal{J}$, then the approximation at this step is $\tilde{A}_z = \tilde{C}_z\tilde{U}_z\tilde{R}_z$, where $\tilde{C}_z = A_{:,\mathcal{J}}$, $\tilde{U}_z = A_{\mathcal{I},\mathcal{J}}^{-1}$, and $\tilde{R}_z = A_{\mathcal{I},:}$, and we need this approximation to search for the next pivot. In most cases, using this directly is fine, but it can suffer from numerical degeneracy from floating point error if $\tilde{U}_z$ is ill-conditioned, especially when $z$ approaches the numerical rank of $A$. Therefore, we derive an iterative construction to avoid the direct action of $\tilde{U}_z$.

We first assume that we already have $z$ row and column indices in $\mathcal{I}$ and $\mathcal{J}$, and a new index pivot $(i,j)$ has been selected such that $A(i,j)\neq 0$. By an index rearrangement, we have the new component matrices in the block form:
\begin{equation*}
    \tilde{C}_{z+1} = \begin{bmatrix}
        \tilde{C}_z & A_{:, j}
    \end{bmatrix},\quad\quad 
    \tilde{U}_{z+1} = \begin{bmatrix}
        \tilde{U}_z^{-1} & A_{\mathcal{I},j}\\
        A_{i,\mathcal{J}} & A_{i,j}
    \end{bmatrix}^{-1},\quad\quad
    \tilde{R}_{z+1} = \begin{bmatrix}
        \tilde{R}_z\\
        A_{i,:}
    \end{bmatrix}.
\end{equation*}
First, focusing on $\tilde{U}_{z+1}$, by block matrix inversion we have
\begin{equation*}
    \tilde{U}_{z+1} = \begin{bmatrix}
        (\tilde{U}_z^{-1} - A_{\mathcal{I},j}A_{i,j}^{-1}A_{i,\mathcal{J}})^{-1} & 0\\0 & \delta_z
    \end{bmatrix}\begin{bmatrix}
        I&-A_{\mathcal{I},j}A_{i,j}^{-1}\\-A_{i,\mathcal{J}}\tilde{U}_z&I
    \end{bmatrix},
\end{equation*}
where $\delta_z = (A_{i,j} - A_{i,\mathcal{J}}\tilde{U}_zA_{\mathcal{I},j})^{-1}$. Now, by the Sherman-Morrison-Woodbury formula, the top left entry of the first matrix can we re-written as
$$
\tilde{U}_z + \tilde{U}_zA_{\mathcal{I},j}(A_{i,j} - A_{i,\mathcal{J}}\tilde{U}_zA_{\mathcal{I},j})^{-1}A_{i,\mathcal{J}}\tilde{U}_z = \tilde{U}_z + \delta_z\tilde{U}_zA_{\mathcal{I},j}A_{i,\mathcal{J}}\tilde{U}_z.
$$
Substituting this in to the block matrix formula for $\tilde{U}_{z+1}$ and computing the matrix-matrix product yields
$$
\tilde{U}_{z+1} = \begin{bmatrix}
\tilde{U}_z + \delta_z \tilde{U}_zA_{\mathcal{I},j}A_{i,\mathcal{J}}\tilde{U}_z & - \delta_z\tilde{U}_zA_{\mathcal{I},j}\\
-\delta_z A_{i,\mathcal{J}}\tilde{U}_z & \delta_z
\end{bmatrix}.
$$
In this way, we can obtain the new approximation $\tilde{A}_{z+1}$ via
\begin{align}
&\tilde{A}_{z+1} = \tilde{C}_{z+1}\tilde{U}_{z+1}\tilde{R}_{z+1} \nonumber \\
& = \begin{bmatrix}
\tilde{C}_z(\tilde{U}_z + \delta_z \tilde{U}_zA_{\mathcal{I},j}A_{i,\mathcal{J}}\tilde{U}_z ) - \delta_z A_{:,j}A_{i,\mathcal{J}}\tilde{U}_z, & -\delta_z \tilde{C}_z\tilde{U}_zA_{\mathcal{I},j} + \delta_z A_{:,j}
\end{bmatrix}\begin{bmatrix}
        \tilde{R}_z\\
        A_{i,:}
    \end{bmatrix} \label{eq:aca_approx}\\
& = \tilde{C}_z(\tilde{U}_z + \delta_z \tilde{U}_zA_{\mathcal{I},j}A_{i,\mathcal{J}}\tilde{U}_z )\tilde{R}_z - \delta_z A_{:,j}A_{i,\mathcal{J}}\tilde{U}_z\tilde{R}_z-\delta_z \tilde{C}_z\tilde{U}_zA_{\mathcal{I},j}A_{i,:} + \delta_z A_{:,j}A_{i,:} \label{eq:aca_approx_tmp}
\end{align}
Upon expanding and factoring~\cref{eq:aca_approx_tmp}, we obtain
\begin{align}
\tilde{A}_{z+1} &= \tilde{A}_z + \delta_z(\tilde{C}_z\tilde{U}_zA_{\mathcal{I},j} - A_{:,j})(A_{i,\mathcal{J}}\tilde{U}_z\tilde{R}_z - A_{i,:}) = \tilde{A}_z + \delta_z(A - \tilde{A}_z)_{:,j}(A - \tilde{A}_z)_{i,:} \nonumber \\
&= \tilde{A}_z + \frac{1}{\tilde{E}_z(i,j)}\tilde{E}_z(:,j)\tilde{E}_z(i,:) \label{eq:aca_approx_E},
\end{align}
where $\tilde{E}_z:=\tilde{A}_z - A$ \reva{is the difference between the target matrix $A$ and current approximation $\tilde{A}_z$}. A similar form of~\cref{eq:aca_approx_E} can be found in~\cite{cortinovis2020maximum}, where the construction is led by an LU decomposition of the approximation. In other words, the construction starts with full data of the original matrix. On the contrary, ours begins with a zero initialization of the approximation.

\begin{figure}[t!]
\centering
    \includegraphics[width=0.5\textwidth]{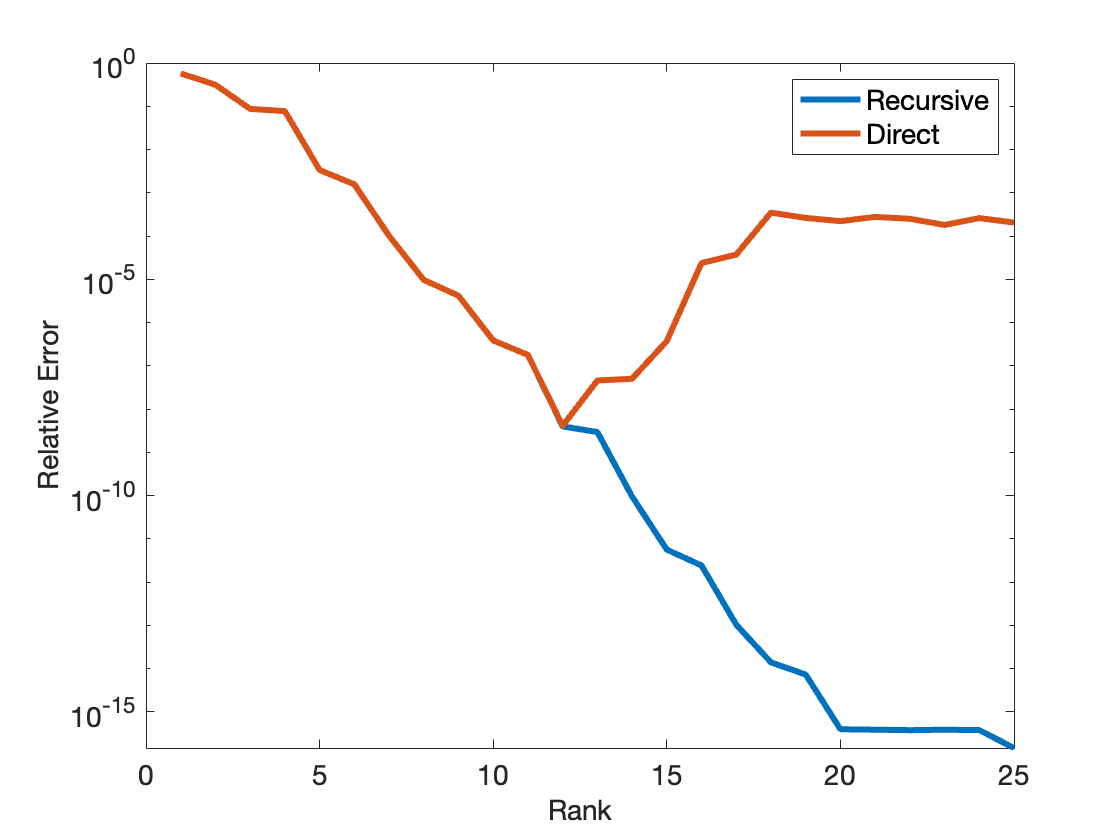}
\caption{Relative Error versus rank using formula \eqref{eq:aca_approx_E} and a direct computation of the factor matrices. The test was run a 100$\times$100 Hilbert matrix.}
\label{fig:Hilbert_recursive_vs_direct}
\end{figure}

As an example, we perform a greedy cross approximation on a 100$\times$100 Hilbert matrix defined elementwise via
\begin{equation} \label{eq:HilbertMatrix}
    H_{i_1,i_2} = \frac{1}{i_1+i_2-1}, \quad 1 \le i_1, i_2 \le 100,
\end{equation}
which is shown to have rapidly decaying singular values in~\cite{beckermann2019bounds}. \reva{In~\cref{fig:Hilbert_recursive_vs_direct}, we compare the behaviors of~\cref{eq:aca_approx_E} and the direct computation for approximation construction, where the evaluation of the inverse of the cross sections $A_{\mathcal{I},\mathcal{J}}$ is used.} There we can see that the use of formula~\cref{eq:aca_approx_E} provides better numerical stability up to machine precision, while the direct method suffers from degeneracy when selected rank is larger than 12. A separate advantage of~\cref{eq:aca_approx_E} is that we are able to make global approximation updates at a local level, which is essential for communication efficiency of submatrix and subtensor algorithms discussed in the following subsections.  

\subsubsection{Iterative construction of submatrix cross approximation 
}
\label{sec:itersubMat1}
Suppose the matrix $A$ is partitioned into $C_1C_2$ submatrices, and we use $C_1C_2$ working processes labeled by $P_{k,\ell}$ for $1 \le k \le C_1$ and $1 \le \ell \le C_2$. This results in a 2D process grid to partition $A$. Using~\cref{eq:aca_approx_E} and the same notations, we can derive an update formula for the submatrix labeled by $(\mathcal{K},\mathcal{L})$ on the process $P_{k,\ell}$
\begin{equation} \label{eq:submatrix_cross_iter}
   \tilde{A}_{z+1}(\mathcal{K},\mathcal{L}) = \tilde{A}_z(\mathcal{K},\mathcal{L})+\frac{1}{\tilde{E}_z(i,j)}\tilde{E}_z(\mathcal{K},j)\tilde{E}_z(i,\mathcal{L}). 
\end{equation} 
This formulation indicates that the construction of the approximation in the new iteration relies on $\tilde{E}_z(i,j)$, $\tilde{E}_z(\mathcal{K},j)$, and $\tilde{E}_z(i,\mathcal{L})$, which might not belong to $P_{k,\ell}$. In this way, there are four main cases:
\begin{enumerate}[leftmargin=*,noitemsep]
\item When $i \in \mathcal{K}$ and $j \in \mathcal{L}$: $\tilde{A}_{z+1}(\mathcal{K},\mathcal{L})$ can be constructed without communications.
\item When $i \in \mathcal{K}$ and $j \notin \mathcal{L}$: $P_{k,\ell}$ needs to obtain $\tilde{E}_z(i,j)$ and $\tilde{E}_z(\mathcal{K},j)$ from $P_{k,\ell^*}$, whose responsible domain is $(\mathcal{K},\mathcal{L}^*)$ with $j \in \mathcal{L}^*$.
\item When $i \notin \mathcal{K}$ and $j \in \mathcal{L}$: $P_{k,\ell}$ needs to obtain $\tilde{E}_z(i,j)$ and $\tilde{E}_z(i,\mathcal{L})$ from $P_{k^*,\ell}$, whose responsible domain is $(\mathcal{K}^*,\mathcal{L})$ with $i \in \mathcal{K}^*$.
\item When $i \notin \mathcal{K}$ and $j \notin \mathcal{L}$: $P_{k,\ell}$ needs to obtain $\tilde{E}_z(i,j)$ from $P_{k^*,\ell^*}$, $\tilde{E}_z(i,\mathcal{L})$ from $P_{k^*,\ell}$, and $\tilde{E}_z(\mathcal{K},j)$ from $P_{k,\ell^*}$, where $k^*$ and $\ell^*$ are responsible for the domain $\mathcal{K}^*$ and $\mathcal{L}^*$ respectively, with $i \in \mathcal{K}^*$ and $j \in \mathcal{L}^*$.
\end{enumerate}

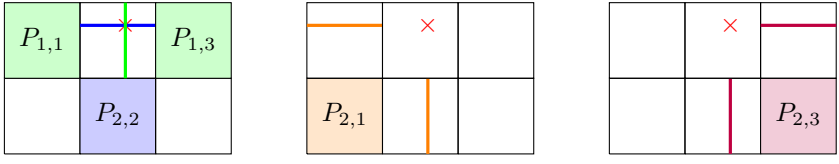
\begin{figure}
    \centering
    \begin{tikzpicture}
        \draw (0,0) -- (1,0) -- (1,1) -- (0,1) -- (0,0);
        \draw[fill=blue!20] (1,0) -- (2,0) -- (2,1) -- (1,1) -- (1,0);
        \filldraw (1.5,0.5) node {$P_{2,2}$};
        \draw (2,0) -- (3,0) -- (3,1) -- (2,1) -- (2,0);
        \draw[fill=green!20] (0,1) -- (1,1) -- (1,2) -- (0,2) -- (0,1);
        \filldraw (0.5,1.5) node {$P_{1,1}$};
        \draw (1,1) -- (2,1) -- (2,2) -- (1,2) -- (1,1);
        \filldraw[red] (1.6,1.7) node {$\times$};
        \draw[blue,very thick] (1,1.7) -- (2,1.7);
        \draw[green,very thick] (1.6,1) -- (1.6,2);
        \draw[fill=green!20] (2,1) -- (3,1) -- (3,2) -- (2,2) -- (2,1);
        \filldraw (2.5,1.5) node {$P_{1,3}$};
    
        \draw[fill=orange!20] (4,0) -- (5,0) -- (5,1) -- (4,1) -- (4,0);
        \filldraw (4.5,0.5) node {$P_{2,1}$};
        \draw (5,0) -- (6,0) -- (6,1) -- (5,1) -- (5,0);
        \draw[orange,very thick] (5.6,0) -- (5.6,1);
        \draw (6,0) -- (7,0) -- (7,1) -- (6,1) -- (6,0);
        \draw (4,1) -- (5,1) -- (5,2) -- (4,2) -- (4,1);
        \draw[orange,very thick] (4,1.7) -- (5,1.7);
        \draw (5,1) -- (6,1) -- (6,2) -- (5,2) -- (5,1);
        \filldraw[red] (5.6,1.7) node {$\times$};
        \draw (6,1) -- (7,1) -- (7,2) -- (6,2) -- (6,1);

        \draw (8,0) -- (9,0) -- (9,1) -- (8,1) -- (8,0);
        \draw (9,0) -- (10,0) -- (10,1) -- (9,1) -- (9,0);
        \draw[purple,very thick] (9.6,0) -- (9.6,1);
        \draw[fill=purple!20] (10,0) -- (11,0) -- (11,1) -- (10,1) -- (10,0);
        \filldraw (10.5,0.5) node {$P_{2,3}$};
        \draw (8,1) -- (9,1) -- (9,2) -- (8,2) -- (8,1);
        \draw (9,1) -- (10,1) -- (10,2) -- (9,2) -- (9,1);
        \filldraw[red] (9.6,1.7) node {$\times$};
        \draw (10,1) -- (11,1) -- (11,2) -- (10,2) -- (10,1);
        \draw[purple,very thick] (10,1.7) -- (11,1.7);
    \end{tikzpicture}
    \caption{The communication pattern for submatrix cross iterative construction if we use a $2 \times 3$ process grid for both matrix $A$ and error $\tilde{E}_z$ at each step. The new pivot is labeled by the red cross on process $P_{1,2}$. In order to compute~\cref{eq:submatrix_cross_iter}, the labeled processes need to receive the highlighted elements (marked with the same color) from neighbour processes.}
    \label{fig:submatrix_cross_comm}
\end{figure}

\Cref{fig:submatrix_cross_comm} illustrates this communication procedure when we use a $2 \times 3$ process grid. In this example, we suppose the new pivot, labeled by the red cross, is found on process $P_{1,2}$, so $P_{1,2}$ belongs to case 1, and $\tilde{A}_{z+1}$ can be constructed with information local to this process. Processes $P_{1,1}$ and $P_{1,3}$ belong to case 2, so they need to receive the column highlighted in green from $P_{1,2}$ (see~\cref{fig:submatrix_cross_comm} (left)). In addition, the situation of $P_{2,2}$ is also depicted in~\cref{fig:submatrix_cross_comm} (left). $P_{2,2}$ is in case 3, so it needs the blue row of $P_{1,2}$ to compute $\tilde{A}_{z+1}$. Finally, processes $P_{2,1}$ (see~\cref{fig:submatrix_cross_comm} (middle)) and $P_{2,3}$ (see~\cref{fig:submatrix_cross_comm} (right)) belong to case 4, so they need to obtain highlighted rows and columns from their neighbor processes respectively.

As a result, building a new submatrix approximation involves at most communications among processes on the same block row and column for two vectors, in addition to one single scalar from the process that handles the new pivot. This distributed version of algorithm is described in~\cref{alg:submatrixCross}, with the assumption that one process handles one submatrix. In practice, our working codes do not need to go through these conditional branches as we can set up sub-communicators for information transfer across processes on the same column or row. Finally, since we use \textit{Allgather} to find the best global pivot, the communication of $\tilde{E}_z(i^*,j^*)$ is thus not needed in~\cref{alg:submatrixCross}. 

\begin{algorithm}
\caption{One step of the matrix cross interpolation algorithm on one process.}
\begin{algorithmic}[1]
\label{alg:submatrixCross}
\Require {Sets $(\mathcal{I},\mathcal{J})$ containing existing pivots, sets $(\mathcal{K},\mathcal{L})$ containing indices of submatrices handled by this process, matrix elements $A(\mathcal{K},\mathcal{L})$, and approximation from the previous step $\tilde{A}(\mathcal{K},\mathcal{L})$.}
\Ensure {A new pivot $(i^*,j^*)$ for the entire matrix $A$.}
\State Perform~\cref{alg:MatCross} on $A(\mathcal{K},\mathcal{L})$ and $\tilde{A}(\mathcal{K},\mathcal{L})$ to get a local pivot $(i^*_p,j^*_p)$.
\State Compute $\tilde{E}(\mathcal{K},\mathcal{L}) = A(\mathcal{K},\mathcal{L}) - \tilde{A}(\mathcal{K},\mathcal{L})$.
\State Use \textit{Allgather} to find the best pivot $(i^*,j^*)$ on all processes.
\If {$i^* \in \mathcal{K}$ and $j^* \in \mathcal{L}$ (Case 1)}
\State \textit{Send} $\tilde{E}(i^*,\mathcal{L})$ to column neighbors (Processes in Case 3).
\State \textit{Send} $\tilde{E}(\mathcal{K},j^*)$ to row neighbors (Processes in Case 2).
\ElsIf {$i^* \in \mathcal{K}$ and $j^* \notin \mathcal{L}$ (Case 2)}
\State \textit{Send} $\tilde{E}(i^*,\mathcal{L})$ to column neighbors (Processes in Case 4).
\State \textit{Receive} $\tilde{E}(\mathcal{K},j^*)$ from a row neighbor that owns the pivot (Case 1).
\ElsIf {$i^* \notin \mathcal{K}$ and $j^* \in \mathcal{L}$ (Case 3)}
\State \textit{Send} $\tilde{E}(\mathcal{K},j^*)$ to row neighbors (Processes in Case 4).
\State \textit{Receive} $\tilde{E}(i^*,\mathcal{L})$ from a column neighbor that owns the pivot (Case 1).
\ElsIf {$i^* \notin \mathcal{K}$ and $j^* \notin \mathcal{L}$ (Case 4)}
\State \textit{Receive} $\tilde{E}(i^*,\mathcal{L})$ from a column neighbor (Process in Case 2).
\State \textit{Receive} $\tilde{E}(\mathcal{K},j^*)$ from a row neighbor (Process in Case 3).
\EndIf
\State Compute the new approximation $\tilde{A}(\mathcal{K},\mathcal{L})$.
\State Set $\mathcal{I} \leftarrow \mathcal{I}\cup i^*$ and $\mathcal{J} \leftarrow \mathcal{J}\cup j^*$.
\end{algorithmic}   
\end{algorithm}

\subsubsection{An alternative iterative construction of matrix and submatrix cross approximation}
\label{sec:iterMat2}
\begin{figure}[h!]
    \centering
    \begin{tikzpicture}[scale=1.25]
        \draw[fill,red,opacity=0.5] (0,0) rectangle (0.2,1.5);
        \node[] at (0.1,-0.25) {\scriptsize{$\tilde{T}_1$}};
        \draw[->] (0.4,0.75)--(0.6,0.75);
        \draw[fill,red,opacity = 0.5] (0.8,0) rectangle (1,1.5);
        \draw[fill,blue,opacity = 0.5] (1,0) rectangle (1.2,1.5);
        \node[] at (1,-0.25) {\scriptsize{$\tilde{T}_2$}};
        \draw[->] (1.4,0.75)--(1.6,0.75);
        \draw[fill,red,opacity = 0.5] (1.8,0) rectangle (2.2,1.5);
        \draw[fill,blue,opacity = 0.5] (2.2,0) rectangle (2.4,1.5);
        \node[] at (2.1,-0.25) {\scriptsize{$\tilde{T}_3$}};
        \draw[->] (2.6,0.75) -- (2.8,0.75);
        \node[] at (3.0,0.75) {...};
        \draw[->] (3.2,0.75)--(3.4,0.75);
        \draw[fill,red,opacity = 0.5] (3.6,0) rectangle (4.6,1.5);
        \draw[fill,blue,opacity = 0.5] (4.6,0) rectangle (4.8,1.5);
        \node[] at (4.2,-0.25) {\scriptsize{$\tilde{T}_{k}$}};
        \draw[->] (5.0,0.75)--(5.2,0.75);
        \draw[fill,red,opacity = 0.5] (5.4,0) rectangle (6.6,1.5);
        \draw[fill,blue,opacity = 0.5] (6.6,0) rectangle (6.8,1.5);
        \node[red] at (6.1,-0.25) {\scriptsize{$\tilde{T}_{k}+\delta_{k}s_{k}\tilde{T}_{k}(i,:)$}};
        \node[blue] at (6.65,1.65) {\scriptsize{$-\delta_ks_k$}};
        \draw[->] (7.0,0.75)--(7.2,0.75);
        \node[] at (7.4,0.75) {...};
        \draw[->] (7.6,0.75)--(7.8,0.75);
        \draw[fill,red,opacity = 0.5] (8,0) rectangle (9.4,1.5);
        \draw[fill,blue,opacity = 0.5] (9.4,0) rectangle (9.6,1.5);
        \node[] at (8.8,-0.25) {\scriptsize{$\tilde{T}_z$}};

    \end{tikzpicture}
    \caption{\reva{Visual representation of the use of~\cref{eq:aca_approx_T} to construct $\tilde{T}_z$. At each iteration a new matrix is constructed by computing the red block using the first term of~\cref{eq:aca_approx_T}, and then the blue column given by the second term of~\cref{eq:aca_approx_T} is concatenated on the right.}}
    \label{fig:Tkconstructiondiagram}
\end{figure}
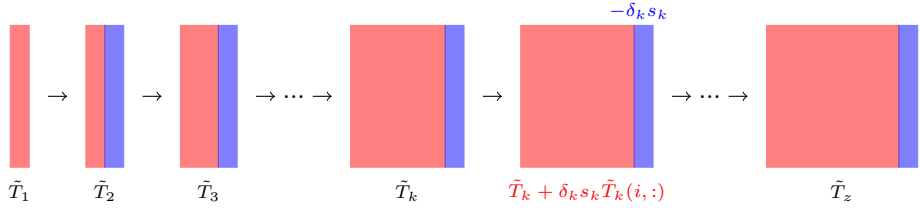

For submatrix ACA,~\cref{eq:submatrix_cross_iter} is sufficient for both finding pivots and constructing approximations. However, building TT cores with only the formulation~\cref{eq:submatrix_cross_iter} is not ideal because the action of the inverse in~\cref{eq:ttcross_build_elem} shall bring the same numerical issues as before into the problem. To overcome this, we derive a recursive formula for $\tilde{T}_z = \tilde{C}_z\tilde{U}_z$, which can be generalized to dimension-wise TT cores in the tensor setting \reva{as they take the form of $CU$. We start by taking the first term in~\cref{eq:aca_approx} and expand and substitute $\tilde{T}_z$ where it appears.
\begin{align*}
\tilde{C}_{z+1}\tilde{U}_{z+1}  & = \begin{bmatrix}
\tilde{C}_z(\tilde{U}_z + \delta_z \tilde{U}_zA_{\mathcal{I},j}A_{i,\mathcal{J}}\tilde{U}_z ) - \delta_z A_{:,j}A_{i,\mathcal{J}}\tilde{U}_z, & -\delta_z \tilde{C}_z\tilde{U}_zA_{\mathcal{I},j} + \delta_z A_{:,j}
\end{bmatrix}\\
& = \begin{bmatrix}
    \tilde{T}_z + (\delta_z \tilde{T}_zA_{\mathcal{I},j}  - \delta_z A_{:,j})A_{i,\mathcal{J}}\tilde{U}_z, & -\delta_z \tilde{T}_zA_{\mathcal{I},j} + \delta_z A_{:,j}
\end{bmatrix}.
\end{align*}
If we then define $s_z= \tilde{T}_zA_{\mathcal{I},j} - A_{:,j}$, and note that $A_{i,\mathcal{J}}\tilde{U}_z = \tilde{C}_z(i,:)\tilde{U}_z = \tilde{T}_z(i,:)$, then we have the recursive formula for $\tilde{T}_{z+1}$:
\begin{equation} \label{eq:aca_approx_T}
    \tilde{T}_{z+1} = \left[\tilde{T}_z + \delta_zs_z\tilde{T}_z(i,:), \quad -\delta_zs_z\right].
\end{equation}
We include~\cref{fig:Tkconstructiondiagram} to demonstrate the functionality of~\cref{eq:aca_approx_T}. In this figure, we see that at each iteration, the red block grows one column at a time, and it depends on the previous iterate via the formula $\tilde{T}_z + \delta_zs_z\tilde{T}_z(i,:)$. This new altered block is then concatenated on the right by a new column vector $-\delta_zs_z$ (shown in blue), and then the full block matrix is used to complete the next step of iteration.}

 In terms of submatrix $(\mathcal{K},\mathcal{L})$, $\mathcal{L}$ contributes to the selection of rows $\tilde{R}_z$, so $\tilde{T}_z$ can be partitioned with respect to only $\mathcal{K}$:
\begin{equation} \label{eq:aca_approx_T_sub}
    \tilde{T}_{z+1}(\mathcal{K},:) = \left[\tilde{T}_z(\mathcal{K},:) + \delta_zs_z(\mathcal{K})\tilde{T}_z(i,:),\quad -\delta_zs_z(\mathcal{K})\right],
\end{equation}
where $s_z(\mathcal{K}) = \tilde{T}_z(\mathcal{K},:)A_{\mathcal{I},j}-A_{\mathcal{K},j}$. \reva{We would like to emphasize that in~\cref{eq:aca_approx_T_sub}, we only need one group of index sets $\mathcal{K}$ to partition the matrices $\tilde{T}_z$ into submatrices that contain all the columns.} This suggests that, instead of using the same 2D grid as in the pivot selection stage as depicted in \cref{fig:submatrix_cross_comm}, a 1D process grid is sufficient for the construction of $\tilde{T}_z$, so it is most convenient to use~\cref{eq:aca_approx_T_sub} once all pivots for the cross approximation are discovered.~\Cref{alg:submatrixCrossnew} describes this procedure of building $T$ 
on a 1D process grid with all pivots selected, and we shall use this algorithm in the TT case for all dimensions. \reva{Simply speaking,~\cref{alg:submatrixCrossnew} iteratively uses~\cref{eq:aca_approx_T_sub} to build a submatrix $T(\mathcal{K},:)$ of the final approximation $T$ in $A = TR$ on a distributed memory computing environment, via plugging the discovered pivots one at a time. Based on whether $i_z$, the row index of the used pivot, is in $\mathcal{K}$, this process either computes $\delta$ or receives $\delta$ and $T(i_z,:)$ from another process.} Furthermore, as the 2D process grid consists of much more processes than the 1D grid, we only need to select a few to build the 1D grid, as opposed to introducing more processes into our algorithm. We can also choose the processes so that communication to gain the information for building $T$ is minimal. 

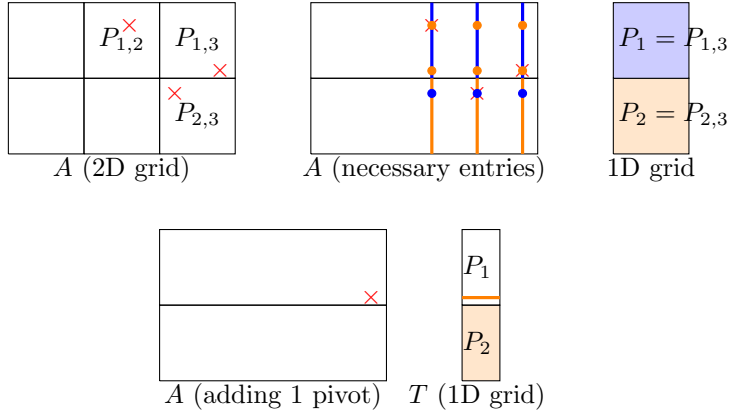
\begin{figure}
    \centering
    \begin{tikzpicture}
        \draw (0,0) -- (1,0) -- (1,1) -- (0,1) -- (0,0);
        \draw (1,0) -- (2,0) -- (2,1) -- (1,1) -- (1,0);
        \draw (2,0) -- (3,0) -- (3,1) -- (2,1) -- (2,0);
        \filldraw[red] (2.2,0.8) node {$\times$};
        \draw (0,1) -- (1,1) -- (1,2) -- (0,2) -- (0,1);
        \draw (1,1) -- (2,1) -- (2,2) -- (1,2) -- (1,1);
        \filldraw[red] (1.6,1.7) node {$\times$};
        \draw (2,1) -- (3,1) -- (3,2) -- (2,2) -- (2,1);
        \filldraw[red] (2.8,1.1) node {$\times$};
        \filldraw[black] (1.5,-0.2) node {$A$ (2D grid)};
        \filldraw (1.5,1.5) node {$P_{1,2}$};
        \filldraw (2.5,1.5) node {$P_{1,3}$};
        \filldraw (2.5,0.5) node {$P_{2,3}$};
    
        \draw (4,0) -- (7,0) -- (7,1) -- (4,1) -- (4,0);
        \draw (4,1) -- (7,1) -- (7,2) -- (4,2) -- (4,1);
        \filldraw[red] (6.2,0.8) node {$\times$};
        \filldraw[red] (5.6,1.7) node {$\times$};
        \filldraw[red] (6.8,1.1) node {$\times$};
        \draw[blue,very thick] (5.6,1) -- (5.6,2);
        \draw[blue,very thick] (6.2,1) -- (6.2,2);
        \draw[blue,very thick] (6.8,1) -- (6.8,2);
        \draw[orange,very thick] (5.6,0) -- (5.6,1);
        \draw[orange,very thick] (6.2,0) -- (6.2,1);
        \draw[orange,very thick] (6.8,0) -- (6.8,1);
        \filldraw[blue] (6.2,0.8) circle (1.5pt);
        \filldraw[blue] (5.6,0.8) circle (1.5pt);
        \filldraw[blue] (6.8,0.8) circle (1.5pt);
        \filldraw[orange] (6.2,1.1) circle (1.5pt);
        \filldraw[orange] (5.6,1.1) circle (1.5pt);
        \filldraw[orange] (6.8,1.1) circle (1.5pt);
        \filldraw[orange] (6.2,1.7) circle (1.5pt);
        \filldraw[orange] (5.6,1.7) circle (1.5pt);
        \filldraw[orange] (6.8,1.7) circle (1.5pt);
        \filldraw[black] (5.5,-0.2) node {$A$ (necessary entries)};

        \draw[fill=orange!20] (8,0) -- (9,0) -- (9,1) -- (8,1) -- (8,0);
        \draw[fill=blue!20] (8,1) -- (9,1) -- (9,2) -- (8,2) -- (8,1);
        \filldraw[black] (8.8,0.5) node {$P_2=P_{2,3}$};
        \filldraw[black] (8.8,1.5) node {$P_1=P_{1,3}$};
        \filldraw[black] (8.5,-0.2) node {1D grid};


        \draw (2,-3) -- (5,-3) -- (5,-2) -- (2,-2) -- (2,-3);
        \draw (2,-2) -- (5,-2) -- (5,-1) -- (2,-1) -- (2,-2);
        \filldraw[red] (4.8,-1.9) node {$\times$};
        \filldraw[black] (3.5,-3.2) node {$A$ (adding 1 pivot)};

        \draw[fill=orange!20] (6,-3) -- (6.5,-3) -- (6.5,-2) -- (6,-2) -- (6,-3);
        \draw (6,-2) -- (6.5,-2) -- (6.5,-1) -- (6,-1) -- (6,-2);
        \filldraw[black] (6.2,-3.2) node {$T$ (1D grid)};
        \filldraw[black] (6.2,-1.5) node {$P_1$};
        \filldraw[black] (6.2,-2.5) node {$P_2$};
        \draw[orange,very thick] (6,-1.9) -- (6.5,-1.9);
    \end{tikzpicture}
    \caption{The communication pattern of using~\cref{alg:submatrixCrossnew} to build approximation~\cref{eq:aca_approx_T_sub}. Top left: three pivots are found with a $2 \times 3$ grid and are denoted by red crosses. Top middle: matrix information (both lines and dots of two colors) needed for the construction of $T$. Top right: A 1D process grid used for building $T$ if row indices of $A$ are partitioned to two pieces. $P_1$ deals with the top half so it needs all blue entries from top middle plot. $P_2$ needs all orange elements from top middle plot. Process $P_{1,3}$ and $P_{2,3}$ from the 2D grid in top left plot are used as $P_1$ and $P_2$ respectively to minimize communication of getting the necessary elements. Bottom left: suppose~\cref{alg:submatrixCrossnew} is performed and there is only one pivot left. Bottom right: at this step, since $P_1$ holds the row index of the pivot, to build the orange portion of $T$, $P_2$ needs to receive the orange row from $P_1$. 
    }
    \label{fig:submatrix_cross_second_comm}
\end{figure}

\cref{fig:submatrix_cross_second_comm} shows a simple illustration of applying~\cref{alg:submatrixCrossnew} to construct the approximation~\cref{eq:aca_approx_T_sub} on the same example in~\cref{fig:submatrix_cross_comm}. Here, we also use a 2D grid of size $2 \times 3$. {Assume} that we find a total of three pivots. The pivots are labeled as red crosses and are found on three processes. Then, as the row indices are partitioned to two pieces, we can use two processes $P_1$ and $P_2$ to construct $T(\mathcal{K},:)$. As the input of~\cref{alg:submatrixCrossnew}, it is required that each process needs to know $A(\mathcal{K}\cup\mathcal{I},\mathcal{J})$. As depicted in the top middle plot in \cref{fig:submatrix_cross_second_comm}, $P_1$ needs to obtain blue lines and dots of $A$ from the 2D process grid; similarly, $P_2$ needs to obtain orange lines and dots. To minimize communication, we choose to use $P_{1,3}$ as $P_1$ and $P_{2,3}$ as $P_2$ since the column indices of two pivots are in the same column partitioning of $A$. The bottom two pictures show an example of the communications in lines 4-8 in~\cref{alg:submatrixCrossnew} that we handle the pivots one by one. Suppose we finish with two pivots and there is only one remaining (see bottom left), we first determine that the row index is in the possession of $P_1$. This indicates that $P_2$ needs to receive the orange row vector from $P_1$ (see bottom right) to build its portion of $T$. 

\begin{algorithm}
\caption{Submatrix cross approximation using~\cref{eq:aca_approx_T_sub} with all pivots.}
\begin{algorithmic}[1]
\label{alg:submatrixCrossnew}
\Require {Sets $(\mathcal{I},\mathcal{J})=((i_1,j_1),\cdots,(i_{N},j_{N}))$ containing all pivots, sets $\mathcal{K}$ containing row indices handled,
and matrix elements $A(\mathcal{K}\cup\mathcal{I},\mathcal{J})$.
} 
\Ensure {The approximation $T(\mathcal{K},:)$ in~\cref{eq:aca_approx_T_sub}.
}
\State Set $T(\mathcal{K},:)=A(i_1,j_1)^{-1}A(\mathcal{K},j_1)$.
\State Set $\mathcal{I}_z = [i_1]$ and $\mathcal{J}_z = [j_1]$.
\For {$2 \le z \le N$}
\If {$i_z \in \mathcal{K}$}
\State Compute $\delta$ via $\delta^{-1} = A(i_z,j_z)-T(i_z,:)A(\mathcal{I}_z,j_z)$. 
\State \textit{Send} $\delta$ and $T(i_z,:)$ 
\Else
\State \textit{Receive} $\delta$ and $T(i_z,:)$.
\EndIf
\State Compute $s(\mathcal{K}) = T(\mathcal{K},:)A(\mathcal{I}_z,j_z)-A(\mathcal{K},j_z)$.
\State Construct $T(\mathcal{K},:) = [T(\mathcal{K},:)+\delta s(\mathcal{K})T(i_z,:) \quad -\!\delta s(\mathcal{K})]$.
\State Set $\mathcal{I}_z = [\mathcal{I}_z \quad i_z]$ and $\mathcal{J}_z = [\mathcal{J}_z \quad j_z]$.
\EndFor
\end{algorithmic}   
\end{algorithm}

{Below we analyze the communication complexity for \cref{alg:submatrixCrossnew}.} Suppose we select $N$ pivots, then with some explicit counting, we can find that a process in the 1D grid needs to receive at most one $\delta$ value and a row vector of length \reva{$z-1$} in the $z$th iteration. Therefore, for all $N-1$ iterations, one process receives at most
\begin{equation} \label{eq:alg32recv}
    (N-1)+\reva{\sum_{z=2}^{N}(z-1)} = (N-1)(N/2+1)
\end{equation}
elements if it does not execute lines 5-6 at all for all pivots. In addition,~\cref{eq:alg32recv} also provides an upper bound of the number of elements one process needs to send to every single other process, if it is responsible to compute lines 5-6 for all pivots. In practice, since the pivots are usually scattered among the processes in the 1D grid, the actual throughput of~\cref{alg:submatrixCrossnew} is much smaller than~\cref{eq:alg32recv}. Nevertheless,~\cref{eq:alg32recv} is useful to analyze communication patterns of our distributed parallel algorithms, and we shall use the bound again to count the throughput of subtensor TT cross.

\subsection{Subtensor TT cross approximation with iterative construction}
\label{sec:iterTT}
In this subsection, we discuss the subtensor TT cross approximations, built upon the submatrix cross algorithm discussed in the previous subsection. Submatrix cross is helpful to subtensor TT cross, since superblocks in line 2 of~\cref{alg:TTcross} are submatrices of tensor unfoldings. When flattened, subtensors correspond to submatrices of tensor unfoldings of various patterns.

\subsubsection{Pivot selection}
\label{sec:TTcrossPiv}
Since cross approximation for the TT format requires pivots regarding multiple dimensions, we aim to discover pivots with the nestedness property~\cite{dolgov2020parallel}, so that the approximation generated by~\cref{alg:TTcross} recovers the exact same elements of the original tensor on positions $(\mathcal{I}_{\le k},\mathbb{I}_{k+1},\mathcal{J}_{>k+1})$ throughout all iterations. \reva{For a more detailed discussion of how nestedness preserves interpolation and why the greedy pivot selection strategy attains nestedness, we refer the readers to~\cite{dolgov2020parallel,savostyanov2014quasioptimality}.} Mathematically, the nestedness of the pivots can be represented as
\begin{equation} \label{eq:nestedpivots}
\mathcal{I}_{\le k+1} \subset \mathcal{I}_{\le k} \times \mathbb{I}_{k+1}, \quad \mathcal{J}_{>k} \subset \mathcal{J}_{>k+1} \times \mathbb{J}_{k+1},
\end{equation}
for $1 \le k \le d-1$. In subtensor parallelism, we \reva{also} hope to maintain this nestedness \reva{by using~\cref{alg:TTcross} on subtensors instead of the entire tensor}.

To see this, for a $d$-dimensional tensor, we first define $(\mathcal{K}_1^{(c)},\dots,\mathcal{K}_d^{(c)})$ to denote the index sets of any subtensor, for $1 \le c \le C$. Then, the sets of pivots on this subtensor can be represented as
\[ \mathcal{L}_{\le 1}^{(c)} = \mathcal{I}_{\le 1} \cap \mathcal{K}_1^{(c)}, \quad \mathcal{L}_{\le 2}^{(c)} = \mathcal{I}_{\le 2} \cap (\mathcal{K}_1^{(c)} \times \mathcal{K}_2^{(c)}), \quad\cdots, \] 
and similarly for $\mathcal{L}_{> 1}^{(c)} = \mathcal{J}_{> 1} \cap \mathcal{K}_1^{(c)}$ and etc.. In addition, it's straightforward to see
\[ \mathcal{K}_1^{(c)} = \mathbb{I}_1 \cap \mathcal{K}_1^{(c)}, \quad \mathcal{K}_2^{(c)} = \mathbb{I}_2 \cap \mathcal{K}_2^{(c)},\quad \cdots. \]
In this way, the nestedness property on each subtensor can be translated as:
\begin{align*}
\mathcal{L}_{\le k}^{(c)} &= \mathcal{I}_{\le k} \cap \left(\mathcal{K}_1^{(c)} \times\cdots\times \mathcal{K}_k^{(c)}\right) \\
&\subset (\mathcal{I}_{\le k-1} \times \mathbb{I}_k) \cap \left(\mathcal{K}_1^{(c)} \times\cdots\times \mathcal{K}_k^{(c)}\right) \\
&\subset \left[\mathcal{I}_{\le k-1} \cap \left(\mathcal{K}_1^{(c)} \times\cdots\times \mathcal{K}_{k-1}^{(c)}\right)\right] \times \left( \mathbb{I}_k \cap \mathcal{K}_k^{(c)}\right) \\
&= \mathcal{L}_{\le k-1}^{(c)} \times \mathcal{K}_k^{(c)}.
\end{align*}
Using the similar argument, we can show $\mathcal{L}_{>k}^{(c)} \subset \mathcal{L}_{>k+1}^{(c)} \times \mathcal{L}_{k+1}^{(c)}$ as well. These provide nestedness guarantee of the subtensor parallel algorithm. As a major corollary, the global best pivots can be thus obtained as the best of the local pivots on each subtensor. In practice, one could perform the pivot selection in a dimensional parallel manner. After pivots are selected for all $d-1$ tensor unfoldings, we are ready to carry out our next task, i.e. to build the TT cores~\cref{eq:submatrix_cross_iter} from the chosen pivots. 

\subsubsection{TT core construction}
\label{sec:TTcrossCore}
In order to construct the TT cores defined in \eqref{eq:ttcross_build}, we use the formula in \eqref{eq:aca_approx_T_sub} along with array slicing to reduce computational requirements. For notational simplicity, we assume that we deal with a $d$-dimensional tensor $\mathcal{X}$ with uniform mode size $n$ and uniform TT rank $r$, and we define $T_{\le k} \in \R^{n^k \times r}$ to be the approximation constructed if~\cref{alg:submatrixCrossnew} is applied for the $k$th unfolding $X_k$. In this way, the first and last TT cores require no special treatment as formation of the first TT core is simply reshaping $T_{\le 1}$ to the correct dimensions, while the last TT core is just a reshape of the selected rows in the last unfolding $X_{d-1}$.
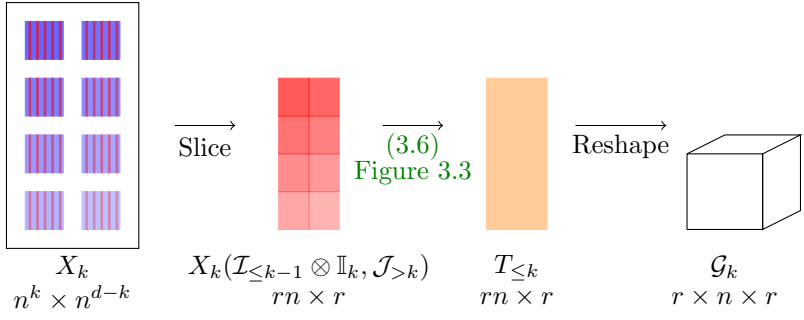
\begin{figure}
    \centering
    \begin{tikzpicture}
        \draw[] (0,-0.25) rectangle (1.75,3);
        \node[] at (0.875,-0.5) {$X_k$};
        \node[] at (0.875,-0.9) {$n^k\times n^{d-k}$};
        \draw[fill,blue,opacity=0.6] (0.25,2.25) rectangle (0.75,2.75);
        \draw[fill,blue,opacity=0.55] (1,2.25) rectangle (1.5,2.75);
        \draw[fill,blue,opacity=0.5] (0.25,1.5) rectangle (0.75,2);
        \draw[fill,blue,opacity=0.45] (1,1.5) rectangle (1.5,2);
        \draw[fill,blue,opacity=0.4] (0.25,0.75) rectangle (0.75,1.25);
        \draw[fill,blue,opacity=0.35] (1,0.75) rectangle (1.5,1.25);
        \draw[fill,blue,opacity=0.3] (0.25,0) rectangle (0.75,0.5);
        \draw[fill,blue,opacity=0.25] (1,0) rectangle (1.5,0.5);
        \draw[->] (2.225,1.375)--(3.025,1.375);

        \draw[fill,red,opacity=0.65] (0.3,2.25) rectangle (0.32,2.75);
        \draw[fill,red,opacity=0.65] (0.4,2.25) rectangle (0.42,2.75);
        \draw[fill,red,opacity=0.65] (0.5,2.25) rectangle (0.52,2.75);
        \draw[fill,red,opacity=0.65] (0.6,2.25) rectangle (0.62,2.75);
        \draw[fill,red,opacity=0.65] (0.7,2.25) rectangle (0.72,2.75);

        \draw[fill,red,opacity=0.6] (0.3+0.75,2.25) rectangle (0.32+0.75,2.75);
        \draw[fill,red,opacity=0.6] (0.4+0.75,2.25) rectangle (0.42+0.75,2.75);
        \draw[fill,red,opacity=0.6] (0.5+0.75,2.25) rectangle (0.52+0.75,2.75);
        \draw[fill,red,opacity=0.6] (0.6+0.75,2.25) rectangle (0.62+0.75,2.75);
        \draw[fill,red,opacity=0.6] (0.7+0.75,2.25) rectangle (0.72+0.75,2.75);

        \draw[fill,red,opacity=0.55] (0.3,1.5) rectangle (0.32,2);
        \draw[fill,red,opacity=0.55] (0.4,1.5) rectangle (0.42,2);
        \draw[fill,red,opacity=0.55] (0.5,1.5) rectangle (0.52,2);
        \draw[fill,red,opacity=0.55] (0.6,1.5) rectangle (0.62,2);
        \draw[fill,red,opacity=0.55] (0.7,1.5) rectangle (0.72,2);

        \draw[fill,red,opacity=0.5] (0.3+0.75,1.5) rectangle (0.32+0.75,2);
        \draw[fill,red,opacity=0.5] (0.4+0.75,1.5) rectangle (0.42+0.75,2);
        \draw[fill,red,opacity=0.5] (0.5+0.75,1.5) rectangle (0.52+0.75,2);
        \draw[fill,red,opacity=0.5] (0.6+0.75,1.5) rectangle (0.62+0.75,2);
        \draw[fill,red,opacity=0.5] (0.7+0.75,1.5) rectangle (0.72+0.75,2);

        \draw[fill,red,opacity=0.45] (0.3,0.75) rectangle (0.32,1.25);
        \draw[fill,red,opacity=0.45] (0.4,0.75) rectangle (0.42,1.25);
        \draw[fill,red,opacity=0.45] (0.5,0.75) rectangle (0.52,1.25);
        \draw[fill,red,opacity=0.45] (0.6,0.75) rectangle (0.62,1.25);
        \draw[fill,red,opacity=0.45] (0.7,0.75) rectangle (0.72,1.25);

        \draw[fill,red,opacity=0.4] (0.3+0.75,0.75) rectangle (0.32+0.75,1.25);
        \draw[fill,red,opacity=0.4] (0.4+0.75,0.75) rectangle (0.42+0.75,1.25);
        \draw[fill,red,opacity=0.4] (0.5+0.75,0.75) rectangle (0.52+0.75,1.25);
        \draw[fill,red,opacity=0.4] (0.6+0.75,0.75) rectangle (0.62+0.75,1.25);
        \draw[fill,red,opacity=0.4] (0.7+0.75,0.75) rectangle (0.72+0.75,1.25);

        \draw[fill,red,opacity=0.35] (0.3,0) rectangle (0.32,0.5);
        \draw[fill,red,opacity=0.35] (0.4,0) rectangle (0.42,0.5);
        \draw[fill,red,opacity=0.35] (0.5,0) rectangle (0.52,0.5);
        \draw[fill,red,opacity=0.35] (0.6,0) rectangle (0.62,0.5);
        \draw[fill,red,opacity=0.35] (0.7,0) rectangle (0.72,0.5);

        \draw[fill,red,opacity=0.3] (0.3+0.75,0) rectangle (0.32+0.75,0.5);
        \draw[fill,red,opacity=0.3] (0.4+0.75,0) rectangle (0.42+0.75,0.5);
        \draw[fill,red,opacity=0.3] (0.5+0.75,0) rectangle (0.52+0.75,0.5);
        \draw[fill,red,opacity=0.3] (0.6+0.75,0) rectangle (0.62+0.75,0.5);
        \draw[fill,red,opacity=0.3] (0.7+0.75,0) rectangle (0.72+0.75,0.5);

        \node[] at (2.625,1.1) {Slice};
        
        \draw[fill,red,opacity=0.65] (3.6,1.5) rectangle (4,2);
        \draw[fill,red,opacity=0.6] (4,1.5) rectangle (4.4,2);
        \draw[fill,red,opacity=0.55] (3.6,1) rectangle (4,1.5);
        \draw[fill,red,opacity=0.5] (4,1) rectangle (4.4,1.5);
        \draw[fill,red,opacity=0.45] (3.6,0.5) rectangle (4,1);
        \draw[fill,red,opacity=0.4] (4,0.5) rectangle (4.4,1);
        \draw[fill,red,opacity=0.35] (3.6,0) rectangle (4,0.5);
        \draw[fill,red,opacity=0.3] (4,0) rectangle (4.4,0.5);
        \node[] at (4,-0.5) {$X_k(\mathcal{I}_{\le k-1}\otimes \mathbb{I}_k,\mathcal{J}_{>k})$};
        \node[] at (4,-0.9) {$rn\times r$};
        \draw[->] (4.975,1.375)--(5.775,1.375);
        \node[] at (5.375,1.1) {\cref{eq:aca_approx_T}};
        \node[] at (5.375,0.75) {\cref{fig:Tkconstructiondiagram}};

        \draw[fill,orange,opacity=0.4] (6.35,0) rectangle (7.15,2);
        \node[] at (6.75,-0.5) {$T_{\leq k}$};
        \node[] at (6.75,-0.9) {$rn\times r$};

        \draw[->] (7.525,1.375)--(8.725,1.375);
        \node[] at (8.125,1.1) {Reshape};

        \draw[] (9,0)--(10,0)--(10,1)--(9,1)--(9,0);
        \draw[] (9,1)--(9.5,1.25)--(10.5,1.25)--(10,1)--(10,0)--(10.5,0.25)--(10.5,1.25);
        \node[] at (9.5,-0.5) {$\mathcal{G}_k$};
        \node[] at (9.5,-0.9) {$r\times n\times r$};
    \end{tikzpicture}
    \caption{The general procedure of constructing an internal core at a global level. \reva{The first image depicts the full superblock $X_k(\mathcal{I}_{\leq k-1}\otimes \mathbb{I}_k,\mathbb{I}_{k+1}\otimes \mathcal{J}_{>k+1})$ for $X_k$ in the blue submatrices. In this image we also show the specific columns $\mathcal{J}_{>k}$ in red, which are sliced out of the superblock to get the second image depicting $X_k(\mathcal{I}_{\leq k-1}\otimes \mathbb{I}_k,\mathcal{J}_{>k})$. This data is given to~\cref{eq:aca_approx_T} (see~\cref{fig:Tkconstructiondiagram}) to obtain $T_{\leq k}$. Finally, $T_{\leq k}$ is reshaped to $\mathcal{G}_k$.}
    }
    \label{fig:core_construction}
\end{figure}
The formation of the internal cores $\mathcal{G}_k$ for $2\leq k\leq d-1$ can start to suffer from computational expense, if we first form the full $T_{\le k}$ and then extract the correct rows determined by the pivots in $\mathcal{I}_{\leq k-1}$. Instead, we implement the reverse by only computing what would be the extracted rows of $T_{\le k}$. In other words, we effectively slice prior to the construction of $T_{\le k}$ so that we only compute a submatrix of $T_{\le k}$. This can be phrased as extraction of a submatrix which contributes to $T_{\le k}$ from the full unfolding $X_k$. Since elements of $\mathcal{I}_{\le k-1}$ populate the first $k-1$ indices of $\mathcal{X}$, we can make this extraction by computing $\mathcal{X}(\mathcal{I}_{\le k-1},:)$. In practice, this can be done efficiently with the function np.ix\_ in numpy by calling $\mathcal{X}[\text{np.ix}\_(\mathcal{I}_{\le k-1})]$.
This reduces the size of the matrix used for computation of $T_{\le k}$ from size $n^{k}\times r$ to size $rn\times r$. In this setting, the neighboring indices $\mathcal{I}_{\leq k-1}$ must be communicated to the processes in charge of the computation of $\mathcal{G}_k$ before the start of construction. Once all processes that require index information for slicing have the necessary information, then we can use~\cref{alg:submatrixCrossnew} to compute the necessary submatrix of $T_{\le k}$, which can then be gathered and reshaped into the global TT cores $\mathcal{G}_k$ using the following relation
\begin{equation}\label{eq:corereshape}
    \mathcal{G}_k = \begin{cases} 
    \text{reshape}(T_{\leq 1},[1,n,r]) & k=1\\
    \text{reshape}(T_{\leq k},[r,n,r]) & 2\leq k\leq d-1\\ \text{reshape}(\mathcal{X}(\mathcal{I}_{d-1},:),[r,n,1])& k=d\end{cases}.
\end{equation}
In summary, the general steps of this process are pictured in~\cref{fig:core_construction}. We extract the necessary entries of the full unfolding $X_k$, which is seen as the transition of the first figure to the second. The second figure is treated with grid type (b) in~\cref{sec:subComm}, and one can find a 3D illustration in~\cref{fig:twogrids_TT}. Then we use~\cref{alg:submatrixCrossnew} to transition from second to third figure. Lastly from third to fourth figure we use~\cref{eq:corereshape} to obtain the TT cores.

\subsubsection{Grid development and overall algorithm}
\label{sec:TTcrossAlg}
With the subroutines of pivot selection and core construction of TT cross, we can assemble the overall algorithm for subtensor TT cross with iterative formulations.

For a $d$-dimensional tensor, we first use~\cref{alg:submatrixCross} with a $d$-dimensional process grid. The processes search for new pivots among the elements in the current superblocks, but build approximations for all the entries in the subtensor. In this way, the approximations of the elements appearing in the superblocks are always updated. This results in more computation to build approximations, but minimizes communications and maintains good load balance in this procedure. Furthermore, since the elements in one subtensor appear in all tensor unfoldings, we can select pivots for all dimensions simultaneously in a dimension parallel manner. Once all the pivots are found, all processes enrich their knowledge of the tensor elements they handle. The overall throughput of this communication is the same as if we enrich every iteration, but with greatly reduced latency.

Then, we follow the discussion in~\cref{sec:TTcrossCore} and use~\cref{alg:submatrixCrossnew} to build the TT cores. In particular, for dimension $k$, $T_{\le k}$ from~\cref{sec:TTcrossCore} is associated with the $k$th unfolding of the tensor, whose row indices represent the combined indices of the first $k$ dimensions of the tensor. Therefore, if each dimension is partitioned to $M$ pieces by the subtensor grid, then we need to build a 1D grid with $M^k$ partitions for $T_{\le k}$. This 1D grid can be simply reshaped and understood as a $k$-dimensional grid (with $M$ pieces per dimension) as a higher-dimensional analogue of the 1D grid used in the matrix case in~\cref{sec:iterMat2}. In practice, we require that one process only appears at most once in any $k$-dimensional grid for $1 \le k \le d-1$ for load balance, and the processes are chosen so that communications to switch between the grids for~\cref{alg:submatrixCross} and~\cref{alg:submatrixCrossnew} are minimal.

\begin{figure}
\centering
\begin{tikzpicture}

\pgfmathsetmacro{\cubex}{2}
\pgfmathsetmacro{\cubey}{2}
\pgfmathsetmacro{\cubez}{2}
\draw[black,fill=gray!20] (0,0,0) -- ++(-\cubex,0,0) -- ++(0,-\cubey,0) -- ++(\cubex,0,0) -- cycle;
\draw[black,fill=gray!20] (0,0,0) -- ++(0,0,-\cubez) -- ++(0,-\cubey,0) -- ++(0,0,\cubez) -- cycle;
\draw[black,fill=gray!20] (0,0,0) -- ++(-\cubex,0,0) -- ++(0,0,-\cubez) -- ++(\cubex,0,0) -- cycle;

\draw[red] (-2,-1) -- (0,-1);
\draw[red] (0,-1) -- (0.77,-0.25);
\draw[red] (-1,0) -- (-1,-2);
\draw[red] (-1,0) -- (-0.3,0.77);
\draw[red] (-1.61,0.39) -- (0.39,0.39);
\draw[red] (0.39,0.39) -- (0.39,-1.61);

\filldraw[black] (-1,-1)  node {$\mathcal{X}$};
\filldraw[black] (-2.3,-1)  node {$\mathcal{K}_1$};
\filldraw[black] (-1,-2.3)  node {$\mathcal{K}_2$};
\filldraw[black] (0.59,-2.01)  node {$\mathcal{K}_3$};

\filldraw[blue] (-1.5,-0.5)  node {$P_1$};
\filldraw[blue] (-1.5,-1.5)  node {$P_2$};
\filldraw[blue] (-0.5,-0.5)  node {$P_3$};
\filldraw[blue] (-0.5,-1.5)  node {$P_4$};
\filldraw[blue] (-1.05,0.55)  node {$P_5$};
\filldraw[blue] (-0.05,0.55)  node {$P_7$};
\filldraw[blue] (0.6,-0.95)  node {$P_8$};

\filldraw[color=black, fill=gray!20] (2,0.8) rectangle (4.5,-1.2);
\filldraw[black] (3.2,-0.2) node {$T_{\le 1}$};
\draw[red] (2,-0.2) -- (4.5,-0.2);
\filldraw[black] (1.7,-0.2)  node {$\mathcal{K}_1$};
\filldraw[blue] (3.2,-0.7)  node {$P_2$};
\filldraw[blue] (3.2,0.3)  node {$P_1$};

\filldraw[color=black, fill=gray!20] (6,0.8) rectangle (7.5,-3.2);
\filldraw[black] (6.7,-1.2) node {$T_{\le 2}$};
\draw[red] (6,-0.2) -- (7.5,-0.2);
\draw[red] (6,-1.2) -- (7.5,-1.2);
\draw[red] (6,-2.2) -- (7.5,-2.2);
\filldraw[black] (8.3,-1.2)  node {$\mathcal{K}_1$\&$\mathcal{K}_2$};
\filldraw[blue] (6.7,0.3)  node {$P_3$};
\filldraw[blue] (6.7,-0.7)  node {$P_4$};
\filldraw[blue] (6.7,-1.7)  node {$P_5$};
\filldraw[blue] (6.7,-2.7)  node {$P_6$};

\end{tikzpicture}
\caption{Grids used for subtensor TT cross on a 3D tensor $\mathcal{X}$. The three dimensions have index partitions $\mathcal{K}_1, \mathcal{K}_2$, and $\mathcal{K}_3$ respectively. Left: $2 \times 2 \times 2$ 3D process and subtensor grid used for pivot selection. Middle: $2 \times 1$ 1D process grid used for core construction of $T_{\le 1}$ associated with the first unfolding $X_1$. Right: $4 \times 1$ 1D process grid used for core construction of $T_{\le 2}$ associated with the second unfolding $X_2$. Since the grid corresponds to two dimensions, it can be considered as a $2 \times 2$ 2D grid. For better load balance, we require no overlap between processes used for $T_{\le 1}$ and those for $T_{\le 2}$. The column size of both $T_{\le 1}$ and $T_{\le 2}$ depends on the iteration number $z$ in~\cref{sec:iterMat2}. 
}
\label{fig:twogrids_TT}
\end{figure}
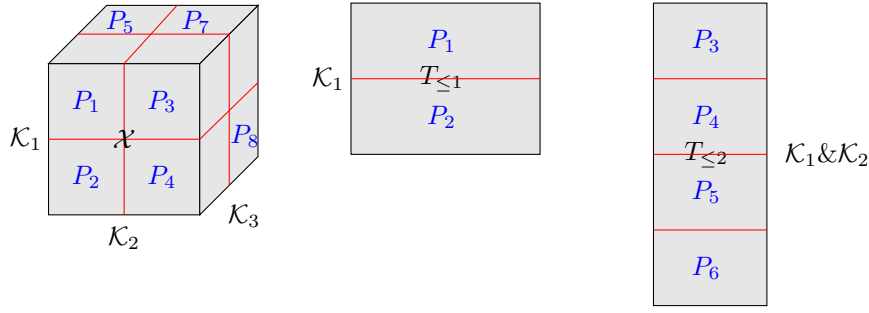

\Cref{fig:twogrids_TT} shows a simple illustration of the three process grids we use for a 3D tensor $\mathcal{X}$ to complete subtensor TT cross. We denote the index partitions of the three dimensions to be $\mathcal{K}_1, \mathcal{K}_2$, and $\mathcal{K}_3$ respectively, and use $P_{p}$ with $1 \le p \le 8$ to represent the active processes. The left plot shows a 3D grid of size $2 \times 2 \times 2$ for both the process and subtensor grid. We use this grid for pivot selection. The middle and right plot show how we partition the construction of the TT cores onto different processes. In particular, $T_{\le 1}$ corresponds to the first unfolding $X_1$ so it needs a $2 \times 1$ 1D grid with respect to $\mathcal{K}_1$. Two processes $P_1$ and $P_2$ partition the construction to two halves. $T_{\le 2}$ corresponds to the second unfolding $X_2$ and thus it needs a $4 \times 1$ 1D grid for both $\mathcal{K}_1$ and $\mathcal{K}_2$, which can be easily reshaped to a 2D grid of size $2 \times 2$ for easier understanding. In this grid, we choose 4 processes other than the two used already for the construction of $T_{\le 1}$. For all three plots, we partition various tasks with red lines, and $P_p$ is used to label the job it works on. For example, $P_1$ in the left figure handles one subtensor out of eight, and constructs the top half part of $T_{\le 1}$ in the middle figure. Furthermore, we want to comment that the column size of both $T_{\le 1}$ and $T_{\le 2}$ depends on the iteration number $z$ in~\cref{sec:iterMat2}. This highlights again that we only need lower dimensional grids to partition $\mathcal{K}_1$ and $\mathcal{K}_2$. 

Finally, we describe this subtensor TT cross approximation with iterative formulations in~\cref{alg:subtensorTTCross}. \reva{In practice,~\cref{alg:subtensorTTCross} is performed with multiple process grids, and more details about data transfers and communication patterns can be found in the next section.}


\begin{algorithm}
\caption{Subtensor TT cross approximation using iterative formulations on one process $c$.}
\begin{algorithmic}[1]
\label{alg:subtensorTTCross}
\Require {Sets $\left(\mathcal{K}_1^{(c)}, \cdots,\mathcal{K}_d^{(c)}\right)$ containing indices of subtensors handled by this process, and tensor elements $\mathcal{X}\left(\mathcal{K}_1^{(c)}, \cdots,\mathcal{K}_d^{(c)}\right)$.}
\Ensure {A TT approximation for $\mathcal{X}$ in distributed memory. Process $c$ returns portion of the TT core $\mathcal{G}_k^{(c)}$ if $c$ is used for constructing the $k$th TT core.}
\For {$0 \le k \le d-1$}
\State Use~\cref{alg:submatrixCross} to find pivot index sets $(\mathcal{I}_{\le k}, \mathcal{J}_{>k})$.
\EndFor
\For {$1 \le k \le d-1$}
\If {$c$ belongs to the grid for computing cores for dimension $k$}
\State \textit{Gather} to get the full $\mathcal{I}_{\le k}$ to the root of this grid.
\If {$c$ is the root of this grid}
\State Receive full $\mathcal{I}_{\le k-1}$ from the root of the grid for dimension $k-1$.
\State Send full $\mathcal{I}_{\le k}$ to the root of the grid for dimension $k+1$.
\EndIf
\State Find pivots in $\mathcal{I}_{\le k-1}$ that appears on $c$ from the root of this grid.
\State Communicate to get tensor elements corresponding to the pivots.
\State Use~\cref{alg:submatrixCrossnew} to construct rows of the approximation $T_{\le k}$.
\State Reshape the constructed rows of $T_{\le k}$ to get $\mathcal{G}_k^{(c)}$.
\EndIf
\EndFor
\end{algorithmic}   
\end{algorithm}

\section{Process grid development for subtensor TT cross}
\label{sec:subComm}
Process grids are often used in distributed memory programming for easier analysis of local computation and communication costs. The idea of 1D and 2D process grids are widely adopted in numerical linear algebra, including matrix multiplications and linear system solvers. In this section, we aim to use process grids for analysis of subtensor parallel TT cross.

For simplicity, for a $d$-dimensional tensor $\mathcal{X}$, we assume it has uniform mode size $n$, subtensor partitioning count per dimension $C$, and process count per dimension $P$. Note that unlike in~\cref{sec:subTTcross} that we assume $P = C$, in this section we no longer have this requirement. In general, one may consider process grids partition the distributed memory machines, while subtensor grids partition the data into optimal size for handling and algorithm performance. In practice, subtensor partitioning is often a refinement of process partitioning.
In addition, we also assume the TT core size is $r$ for all dimensions, so the number of pivots we select for each dimension is $r$ as well. Suppose $n = Cm$, then $m$ is the mode size per dimension of the subtensor. If there are not enough processes to fill a $d$-dimensional process grid, it's easy to see that one can apply any $z$-dimensional grid as long as $z < d$. In this case, mode sizes of certain dimensions of the subtensors are multiples of $m$. However, lower dimensional process grids can lead to load imbalance issues when $C$ is not a multiple of $P$, \reva{since the subtensors cannot have uniform mode size. As a result, some processes compute with more subtensors, or subtensors with larger size}. In this section, we assume $C=Pw$, so that each process handles $w^d$ subtensors\reva{, and the number of elements per dimension held by one process is $wm$.}

As remarked in~\cref{sec:iterTT},~\cref{alg:subtensorTTCross} has two sets of process grids for~\cref{alg:submatrixCross} and~\cref{alg:submatrixCrossnew}: (a) $d$-dimensional grid for pivot selection, \reva{where each process in the grid stores subtensors of the target tensor,} and (b) $d-1$ lower dimensional grids for TT core construction, \reva{where each process is only responsible to generate part of one TT core, and thus needs to store the tensor elements corresponding to the used pivots in~\cref{alg:submatrixCrossnew}. As the process is already part of grid (a) and stores some tensor entries, it needs to receive the necessary data from other processes in grid (a).}

The overall communication costs of subtensor parallel TT cross~\cref{alg:subtensorTTCross} can be understood in three parts:

\begin{itemize}[leftmargin=*,noitemsep]
\item \textbf{Grid (a): $d$-dimensional grid for pivots.} In TT pivot selection,~\cref{alg:submatrixCross} is used $d-1$ times, one for each tensor unfolding. In this step, we need to study both the \textit{Allgather} for pivots, and communications with row and column neighbors. Once we find local pivots on all subtensors, we can compare among the $w^d$ subtensors on the same process to get $P^d$ potential new pivots for each dimension, labeled by $(i_{\le j}^{(p)},i_{>j}^{(p)})$. Here, $1 \le j \le d-1$ is used to denote dimensions, and $1 \le p \le P^d$ is used to denote a specific process. Next, \textit{Allgather} allows all processes to know all pivots $(i_{\le j}^{(p)},i_{>j}^{(p)})$, together with the pivot selection metric. This procedure requires each process to send $3(P^d-1)$ elements for all dimensions, and receive a throughput of $3(P^d-1)$ entries. \reva{In this counting, 3 consists of the two indices of the local pivot and the pivot value, and $P^d-1$ is the total number of other processes to communicate to get the global pivot.} In this way, the best global pivot for a dimension $k$ can be found on each process.

We then calculate the communications with row and column neighbors. For dimension $j$, and at iteration $z$, the process that does not contain the pivot, nor is a column or row neighbor of the process with the pivot, needs to receive the most information. To be exact, it receives 
\begin{equation} \label{eq:matcross1recv}
    (wm)^j+(wm)^{d-j}
\end{equation}
elements, where the two terms are the vectors received from column and row neighbors respectively \reva{(see~\cref{fig:submatrix_cross_comm})}. Comparatively, the process with the pivot needs to send the most information to its \reva{$P^j-1$} row and \reva{$P^{d-j}-1$} column neighbors, of an approximate entry count
\begin{equation} \label{eq:matcross1send}
    (wm)^j(P^j-1)+(wm)^{d-j}(P^{d-j}-1).
\end{equation}

\item \textbf{Connecting Grid (a) and Grid (b).} The gap between Grid (a) and Grid (b) is bridged by lines 5-10 in~\cref{alg:subtensorTTCross}. Lines 5-9 communicate about the pivots with grids for adjacent dimensions, and line 10 get necessary tensor entries for construction. Communications across the lower dimensional grids in lines 5-9 is straightforward, as they only occur for the pivot index sets on the grid roots. For dimension $j$, this means sending $r$ values to dimension $j+1$, and receiving $r$ values from dimension $j-1$. The analysis of line 10 is a bit more complicated. In the extreme case that one process does not contain any of the pivots, it needs to obtain tensor elements associated with the subtensors $\mathcal{X}(\mathcal{I}_{\le j},\mathcal{J}_{>j})$ and $\mathcal{X}(\mathcal{K}_1,\cdots,\mathcal{K}_j,\mathcal{J}_{>j})$ for $1 \le j \le d-1$ \reva{, which contains $r^2$ and $r(wm)^j$ elements respectively}. Therefore, a pessimistic bound of the total number of elements received is
\begin{equation} \label{eq:matcross2recv}
    \sum_{j=1}^{d-1} r^2+r(wm)^j = (d-1)r^2+r\sum_{j=1}^{d-1} (wm)^j
\end{equation}
for all dimensions. In the other extreme case that one process contains all the pivots,~\cref{eq:matcross2recv} is an upper bound for the number of elements sent across all processes. Also, at most $\sum_{j=1}^{d-1}P^j$ processes send and receive in Grid (b).

\item \textbf{Grid (b): Lower-dimensional grids for TT construction.}
For a $d$-way tensor, Grid (b) contains $d-1$ lower-dimensional grids for the row indices of tensor unfoldings, i.e. $j$-dimensional grid for the $j$th unfolding for $1 \le j \le d-1$, as the row indices of this unfolding correspond to $j$ dimensions. The communications only happen in line 11 of~\cref{alg:subtensorTTCross}, within the $j$-dimensional grid. Since one process can only appear once in one of the lower dimensional grids, the bound~\cref{eq:alg32recv} can be applied to the analysis as well. As the final TT core is simply a reshape of the computed portion of $T_{\le j}$, this part of the algorithm does not involve a lot of hanging procedures and communications, and is thus very fast in practice.

\end{itemize}



\section{Numerical Examples}
\label{sec:NumericalExamples}
In this section, we show the performance of our parallel TT cross algorithm\footnote{The code for implementation can be found at \href{https://github.com/dhayes95/Cross}{github.com/dhayes95/Cross}.} on the Hilbert tensors (see~\cref{sec:HilbertTensor}) and discretized Maxwellian functions as equilibrium distribution functions in kinetic theory \cite{cercignani1969mathematical,xiong2016high} (see~\cref{sec:Maxwellian}). \reva{For all tests, the tensors used are represented element-wise using function evaluations. This allows for low memory requirement as the construction used does not require full tensor data, so only entries that are required for construction are computed on the fly.} 

In all numerical tests, results are computed on the DARWIN (Delaware Advanced Research Workforce and Innovation Network) system. The standard partition of this system consists of \reva{two 32-core AMD EPYC 7002 2.80 GHz series processors allowing for 64 total cores per compute node. At the current time, DARWIN does not support node-to-node communication using MPI4PY, so a distributed memory setting is simulated by giving each CPU its own memory.} In addition, for all examples, any mention of error measurement corresponds to the relative error on a uniformly generated set of $R$ random indices $i$
$$
||\mathcal{X} - \tilde{\mathcal{X}}||_{F,R}^2 := \frac{\sum^{R}{ (\mathcal{X}(i) - \tilde{\mathcal{X}}(i)})^2}{\sum^{R}{\mathcal{X}(i)^2}}.
$$
This error approximation allows us to visualize the accuracy of our solver without having to go through all the elements of a large dataset.

\subsection{Testing setup}
\label{sec:Testingsetup}
Given a fixed number of processes, there can be multiple combinations of assigning them to all the dimensions of a tensor. In our tests, the process grid partitions used in the strong scaling tests are obtained by testing all options of partitions and selecting the one which has the fastest run time. In contrast, the partitions of weak scaling tests are selected such that for a given number of MPI ranks, the partitioning is spread evenly across all dimensions.~\cref{tab:Hilbert_MPI_Partition} and~\cref{tab:Max_MPI_Partition} show the partitions used, and entries of these tables show the number of uniform subdivisions of each dimension that are used for a certain number of MPI ranks. For example, a partition of $(2,2,2)$ corresponds to the left image of~\cref{fig:twogrids_TT}.

\reva{In all test cases, the core ranks used are selected experimentally with the cardinality of the pivot sets in ACA. While the Hilbert tensor has theoretical justification for the core ranks \cite{shi2021compressibility}, core ranks of the Maxwellian tensors are theoretically unknown and require the experimental selection.} For all results presented, the algorithm is performed for 10 runs, and the minimum times and errors are taken for plots. For the error plots, all core ranks are incrementally increased until they reach the specified core rank used; e.g., in the $4$-d Maxwellian test in~\cref{sec:Maxwellian} with prescribed rank $(1,10,5,20,1)$, the core ranks used for error tests are $$
(1,2,2,2,1),\dots,(1,5,5,5,1),(1,6,5,6,1),\dots,(1,10,5,11,1),\dots,(1,10,5,20,1).$$

\begin{table}[t!]
\begin{center}
\caption{Partitioning per dimension used for strong and weak scaling of Hilbert tensor.}
\label{tab:Hilbert_MPI_Partition}
\begin{tabular}{ |c|c|c|c|c| } 
\hline 
MPI Ranks & Strong $d = 3$ & Weak $d = 3$ & Strong $d = 6$ & Weak $d=6$\\
 \hline
 1 & $(1,1,1)$ & $(1,1,1)$ & $(1,1,1,1,1,1)$ & $(1,1,1,1,1,1)$ \\ 
 2 & $(1,2,1)$ & $(1,2,1)$ & $(1,1,1,2,1,1)$ & $(1,1,1,1,1,2)$\\ 
 4 & $(1,4,1)$ & $(1,2,2)$ & $(1,1,2,2,1,1)$ & $(1,1,1,1,2,2)$\\ 
 8 & $(1,8,1)$ & $(2,2,2)$ & $(1,1,2,2,2,1)$ & $(1,1,1,2,2,2)$\\
 12 & - & $(2,3,2)$& -&-\\
 16& $(1,16,1)$ & -  & $(1,1,4,2,2,1)$ & $(1,1,2,2,2,2)$\\
 18 & -& $(2,3,3)$ & - & -\\
 27 & -& $(3,3,3)$& -&-\\
 32& $(1,32,1)$ & - & $(1,1,4,4,2,1)$ & $(1,2,2,2,2,2)$\\
 36 & - & $(3,4,3)$& -&-\\
 48 & - & $(3,4,4)$& -&-\\
 64& $(1,64,1)$ & $(4,4,4)$& $(1,2,4,4,2,1)$ & $(2,2,2,2,2,2)$ \\
 \hline
\end{tabular}

\end{center}
\end{table}

\begin{table}[t!]
\begin{center}
\caption{Partitioning per dimension used for strong and weak scaling of Maxwellian tensor.}
\label{tab:Max_MPI_Partition}
\begin{tabular}{ |c|c|c|c|c| } 
\hline 
MPI Ranks & Strong $d=4$ & Weak $ d=4$ & Strong $d=6$ & Weak $d = 6$\\
 \hline
 1 & $(1,1,1,1)$ & $(1,1,1,1)$  & $(1,1,1,1,1,1)$ & $(1,1,1,1,1,1)$  \\ 
 2 & $(1,1,2,1)$ & $(2,1,1,1)$ & $(1,1,1,1,2,1)$ & $(1,1,1,1,1,2)$ \\  
 4 & $(1,2,2,1)$ & $(2,1,2,1)$ & $(1,1,1,2,2,1)$ & $(1,1,1,1,2,2)$ \\
 8 & $(1,2,4,1)$ & $(2,2,2,1)$ & $(1,1,2,2,2,1)$ & $(1,1,1,2,2,2)$ \\
 16& $(1,4,4,1)$ & $(2,2,2,2)$  & $(1,2,2,2,2,1)$ & $(1,1,2,2,2,2)$ \\
 24 & -& $(3,2,2,2)$ & -&- \\
 32& $(1,4,8,1)$ & - & - & -\\
 36 & - & $(3,2,3,2)$& $(1,2,4,2,2,1)$ & $(1,2,2,2,2,2)$\\
 48 & - & $(4,2,3,2)$& - & - \\
 64& $(1,8,8,1)$ & $(4,2,4,2)$ & $(1,2,4,4,2,1)$ & $(2,2,2,2,2,2)$ \\
 \hline
\end{tabular}
\end{center}
\end{table}

\subsection{Hilbert tensor}
\label{sec:HilbertTensor}
In this section, we show the accuracy and some scaling results of our subtensor parallel TT cross for the synthetic dataset--Hilbert tensors, which are higher order analogues of the Hilbert matrices~\cref{eq:HilbertMatrix}. A $d$-dimensional Hilbert tensor can be represented element-wise as
\begin{equation} \label{eq:HilbertTensor}
    \mathcal{X}_{i_1,\dots,i_d} = \frac{1}{1-d+i_1+\cdots+i_d}, \quad 1 \le i_j \le n, \quad 1 \le j \le d,
\end{equation}
where we assume a uniform mode size $n$. It is known that these Hilbert tensors can be accurately approximated by a numerically low TT rank tensor, and the TT ranks $r$ can be estimated a priori given $d$ and $n$~\cite{shi2021compressibility}. We shall use this rank approximation as the number of pivots we select in our tests.

\begin{figure}[t!]
    \centering
    \includegraphics[width=0.33\linewidth]{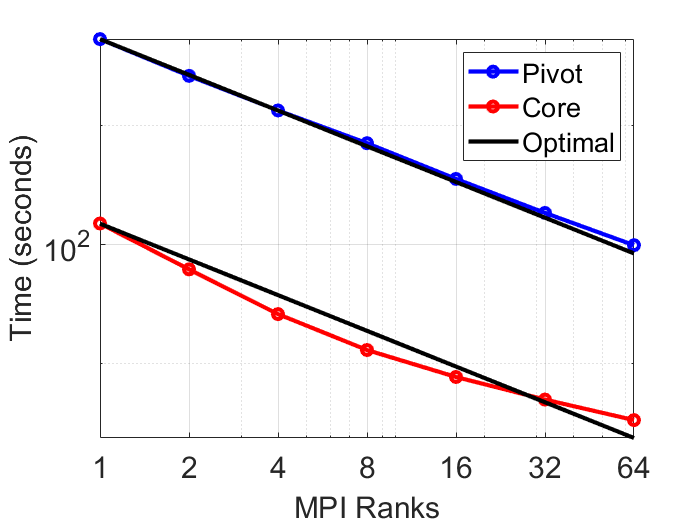}\includegraphics[width=0.33\linewidth]{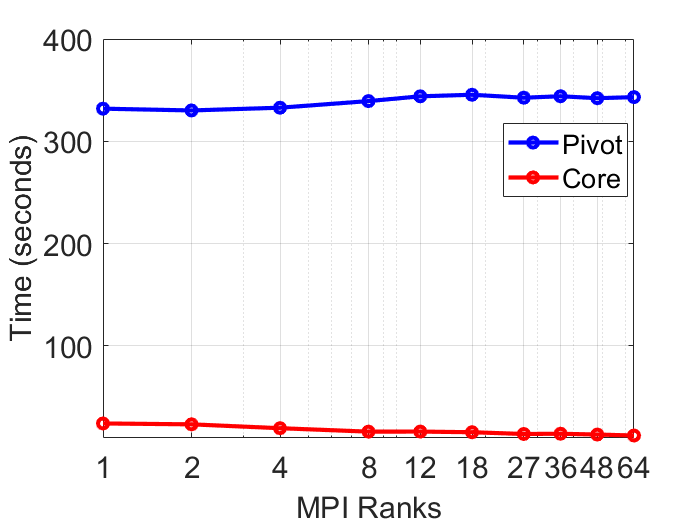}\includegraphics[width=0.33\linewidth]{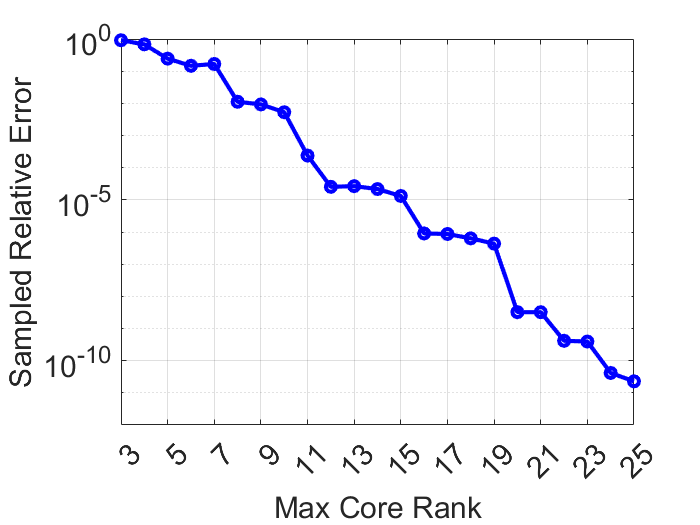}
    \caption{Hilbert tensor with $d = 3$ and $n = 2000$. Left: strong scaling. Middle: weak scaling. Right: accuracy approximation with $||\mathcal{X}- \tilde{\mathcal{X}}||_{F,10^6}$ \reva{using $64$ MPI ranks}.}
    \label{fig:Hilb_3D_Strong_Weak_Err}
\end{figure}

First, we test the case of $d = 3$ with mode size $n = 2000$ and core ranks $(1,25,25,1)$, and refer to~\cref{tab:Hilbert_MPI_Partition} to see partitioning set up for certain MPI ranks in scaling tests. In~\cref{fig:Hilb_3D_Strong_Weak_Err}, we can see the plots for strong and weak scaling results. This test shows good strong scaling results as the pivot selection demonstrates almost optimal results, \reva{where the optimal line corresponds to a line with slope 1, i.e. time is cut in half when the number of MPI ranks is doubled}. For the core construction we see for 1-8 MPI ranks a larger slope. The cause for this phenomena has not currently been investigated, but will be our future target. Following this, we see the scaling from 8-64 MPI ranks demonstrate almost optimal linear scaling. For the weak scaling test, we increase both the mode size and the number of processes linearly, while maintaining a subtensor of size $500\times 500\times 500$ on each process. In~\cref{fig:Hilb_3D_Strong_Weak_Err} we can see in the middle plot that as the number of processes increases we do not suffer from large growth in computational time for pivot selection nor core construction. As with the strong scaling results, the cause behind small intermittent decrease in time is not currently investigated. We also include results of sampled errors in~\cref{fig:Hilb_3D_Strong_Weak_Err} (Right). Here the horizontal axis corresponds to the core ranks $(1,r,r,1)$ and the vertical axis are the values $||\mathcal{X} - \tilde{\mathcal{X}}||_{F,10^6}$. \reva{The core ranks are selected to be $(1,r,r,1)$ as $3D$ Hilbert tensors admit symmetric core ranks from~\cref{eq:HilbertTensor}.} This plot shows that we are able to achieve a close approximation of the true tensor value in the distributed memory setting.

\begin{figure}[t!]
    \centering
    \includegraphics[width=0.5\linewidth]{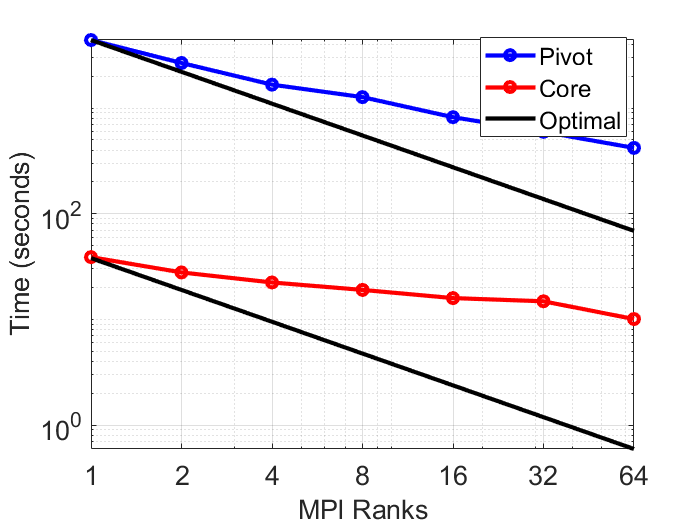}\includegraphics[width=0.5\linewidth]{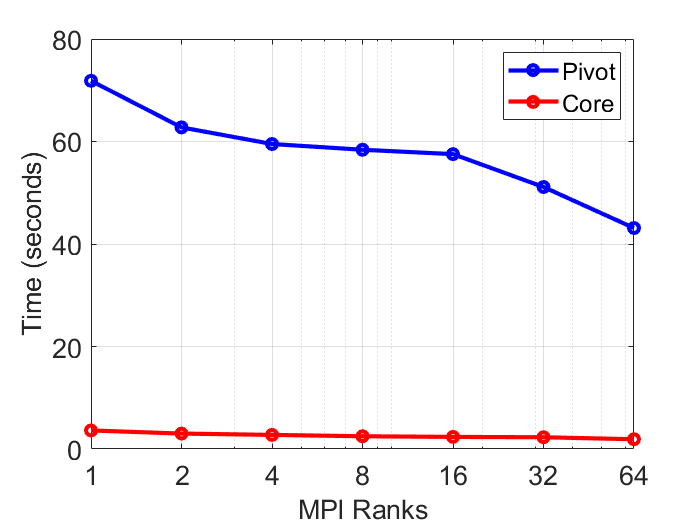}
    \includegraphics[width=0.5\linewidth]{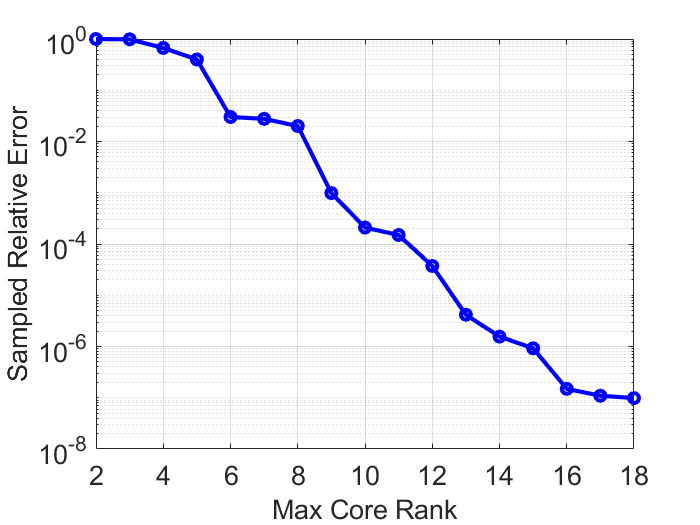}\includegraphics[width=0.5\linewidth]{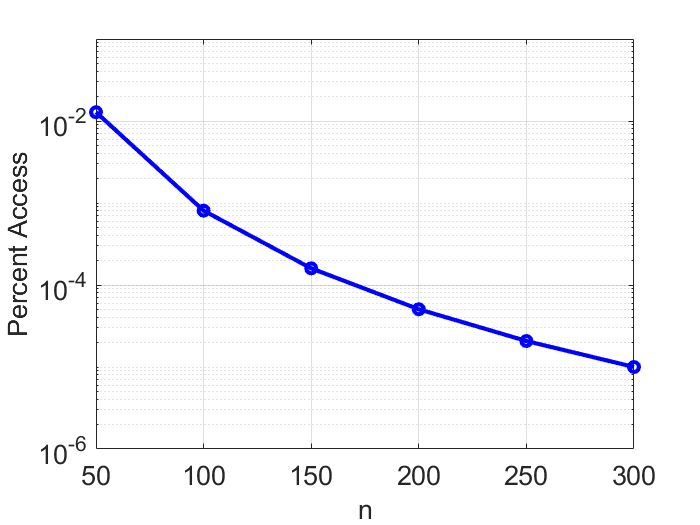}
    \caption{Hilbert tensor with $d = 6$ and $n = 300$. Top left: strong scaling. Top right: weak scaling. Bottom left: accuracy approximation with $||\mathcal{X} - \tilde{\mathcal{X}}||_{F,10^6}$. Bottom right: percentage of full tensor size accessed to perform construction.}
    \label{fig:Hilb_6D_Strong_Weak_Err_Acc}
\end{figure}

For the second test, we again use the Hilbert tensor~\cref{eq:HilbertTensor} with $d = 6$, $n = 300$ and core ranks $(1,15,17,18,17,15,1)$. The MPI partitioning set up can be found in~\cref{tab:Hilbert_MPI_Partition}. Plots for strong and weak scaling are included in~\cref{fig:Hilb_6D_Strong_Weak_Err_Acc} (Top row), and here we see that we can achieve good strong scaling for high dimensional problems. The same phenomena for the $d = 3$ weak scaling results is present for the $d = 6$ case, and also is currently uninvestigated. We suspect that it results from the dimension of the problem handled per caches of the cluster, and the mode size $n=300$ not being large enough. In~\cref{fig:Hilb_6D_Strong_Weak_Err_Acc} (Bottom row) we also include a plots of sampled errors, as well as a measurement of the percentage of the full tensor data that is accessed to construct the TT cores. This shows that in the high dimensional case, we are able to obtain a close approximation. We also can see from the access plot that we require very few true data values, which in turn lends this construction to require very little memory to run. In testing, this local storage requirement ranges from $6$ Gb for 1 MPI rank down to $3.66 \times 10^{-1}$ Gb for 64 MPI ranks to run the full algorithm.

\subsection{Tensors from Maxwellian distribution functions} \label{sec:Maxwellian}
In this section, we show two example tensors, with dimension 4 and dimension 6, constructed from distribution functions arise from kinetic theory of gas dynamics \cite{cercignani1969mathematical,xiong2016high}. In both cases, the tensors are constructed element-wise from function values on a discretized grid in both displacement and velocity spaces. For the $4$-dimensional (2d2v) tensor, the function used to compute entries of $\mathcal{X}$ is given by
\begin{equation} \label{eq:2d2vFunction}
f(x,y,v_x,v_y) = \varrho(x,y)\left[\exp{\left(-b_x^- - b_y^-\right)}+\exp{\left(-b_x^+ - b_y^+\right)}\right]
\end{equation}
where
$$
\varrho(x,y) = \left(\frac{\rho(x)}{2\sqrt{2\pi T(x)}} + \frac{\rho(y)}{2\sqrt{2\pi T(y)}}\right),\quad b_x^{\pm} = \frac{|v_x \pm 0.75|^2}{2T(x)},\quad b_y^{\pm} = \frac{|v_y \pm 0.75|^2}{2T(y)}
$$
and
$$
\rho(w)=1+0.875\sin(2\pi w),\quad T(w)= 0.5 + 0.4\sin(2\pi w), \quad w = x,y.
$$
In our tests, the underlying domains are taken to be $[-1/2,1/2]$ for $x,y$ and $[-3,3]$ for $v_x,v_y$. Then our entries of $\mathcal{X}$ are defined by
$$
\mathcal{X}_{i_1,i_2,i_3,i_4} = f(x_{i_1},y_{i_3},(v_x)_{i_2},(v_y)_{i_4}).
$$

For the $6$-dimensional (3d3v) tensor, the format follows from~\cref{eq:2d2vFunction} with the addition of a $z$ and $v_z$ terms to define $f(x,y,z,v_x,v_y,v_z)$. The specific ordering for the dimensions of $\mathcal{X}$ is $x,v_x,y,v_y$ for $d = 4$, and $x,v_x,y,v_y,z,v_z$ for $d = 6$. These are selected as they provide the best results in practice. Other orderings of the dimensions, e.g. $x,y,v_x,v_y$, are tested, but do not yield small enough errors for large core ranks.

\begin{figure}[t!]
    \centering
    \includegraphics[width=0.33\linewidth]{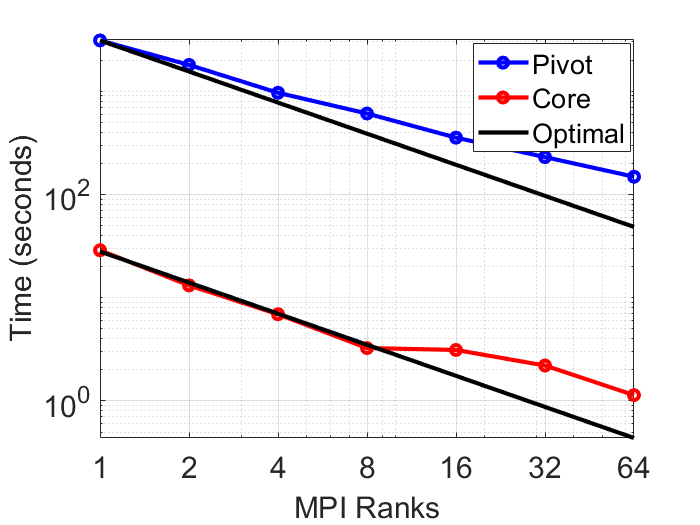}\includegraphics[width=0.33\linewidth]{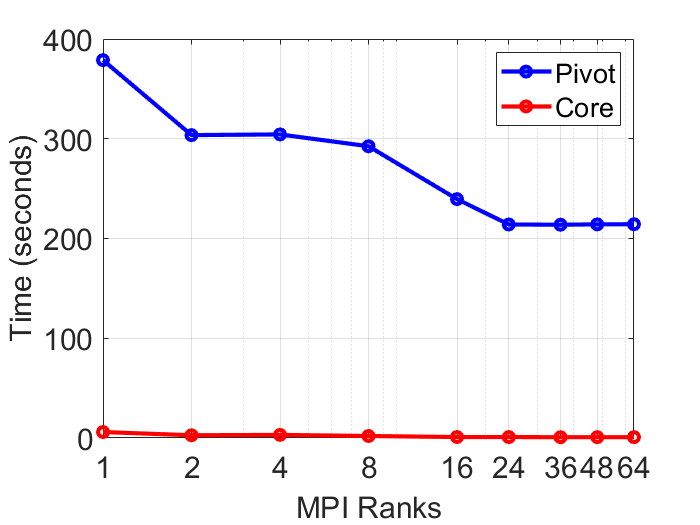}\includegraphics[width=0.33\linewidth]{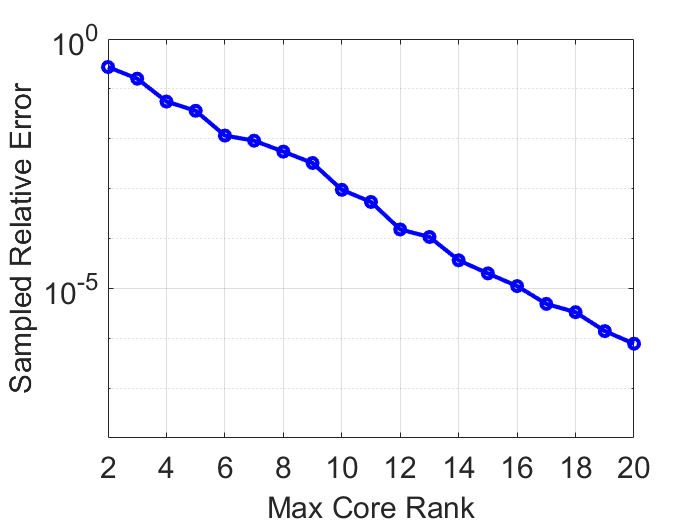}
    \caption{Maxwellian tensor with $d = 4$. Left: strong scaling. Middle: weak scaling. Right: Accuracy approximation with $||\mathcal{X}- \tilde{\mathcal{X}}||_{F,10^6}$ \reva{using 64 MPI ranks.}}
    \label{fig:Max_4D_Strong_Weak_Err}
\end{figure}

For the first test, we work with 2d2v with size $(2n,n,2n,n)$ where $n = 1000$ and core ranks $(1,10,5,20,1)$. As with the Hilbert tensor, the MPI partitioning can be found in~\cref{tab:Max_MPI_Partition}.~\cref{fig:Max_4D_Strong_Weak_Err} shows the plots for both strong and weak scaling for the 2d2v test case. We observe good strong scaling results for all partitions. The weak scaling results display similar results to the Hilbert tensor tests. We also include an error plot of $||\mathcal{X} - \tilde{\mathcal{X}}||_{F,10^6}$, which verifies that we can get a close approximation even with variation in the core ranks as seen in~\cref{fig:Max_4D_Strong_Weak_Err}.

\begin{figure}[t!]
    \centering
    \includegraphics[width=0.5\linewidth]{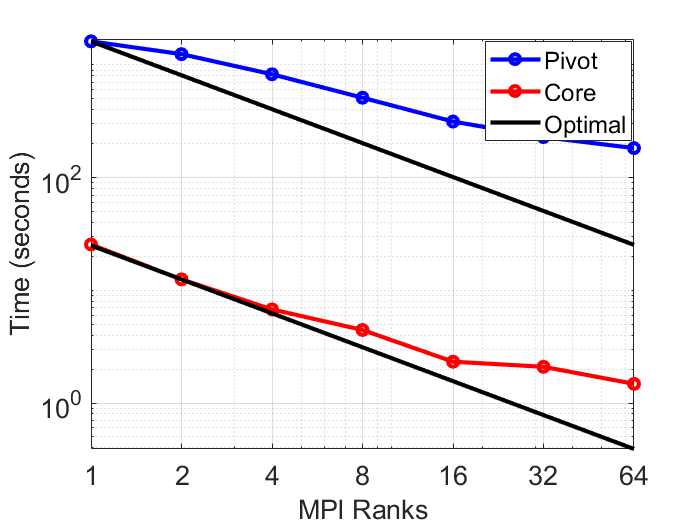}\includegraphics[width=0.5\linewidth]{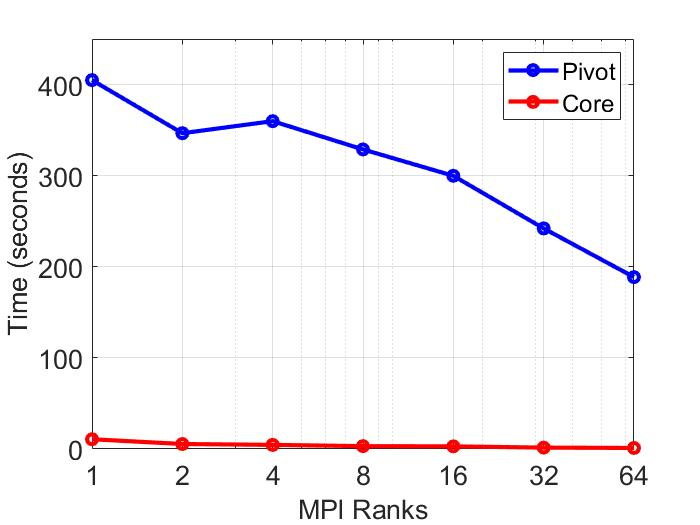}
    \includegraphics[width=0.5\linewidth]{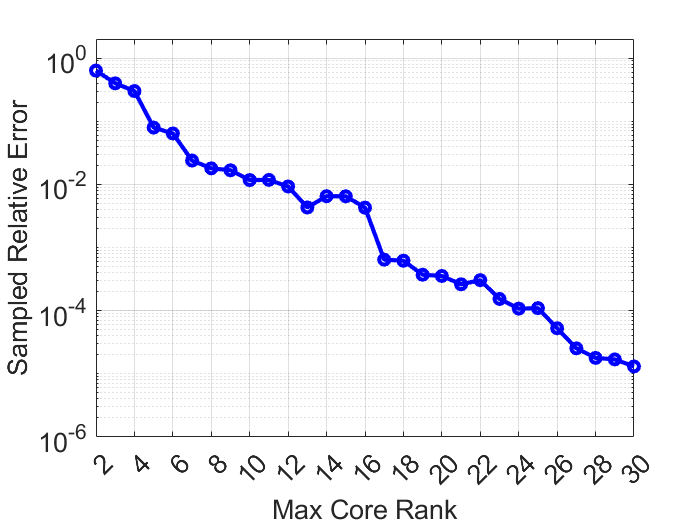}\includegraphics[width=0.5\linewidth]{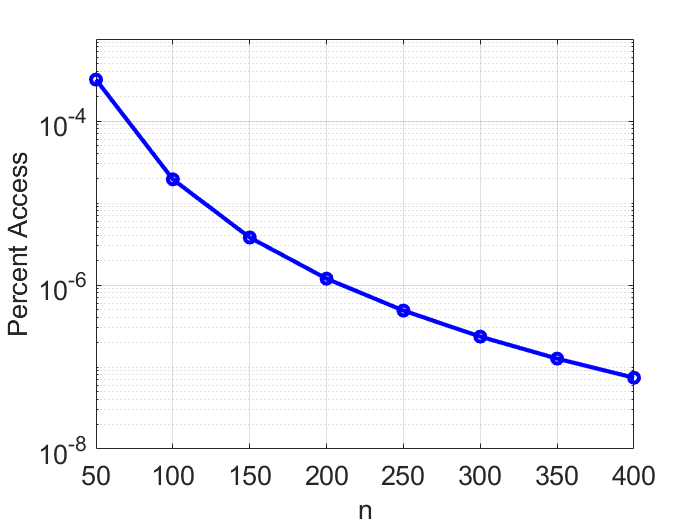}
    \caption{Maxwellian tensor with $d = 6$. Top left: strong scaling. Top right: weak scaling. Bottom left: Accuracy approximation with $||\mathcal{X} - \tilde{\mathcal{X}}||_{F,10^6}$. Bottom right: Percentage of full tensor size accessed to perform construction.}
    \label{fig:Max_6D_Strong_Weak_Err_Acc}
\end{figure}

The last test is on $3d3v$, with size $(2n,n,2n,n,2n,n)$ where $n = 400$, and core ranks $(1,10,5,30,5,20,1)$. The corresponding MPI partitions used are in~\cref{tab:Max_MPI_Partition}. First,~\cref{fig:Max_6D_Strong_Weak_Err_Acc} (Top row) shows the strong and weak scaling results for the given tensor, and we observe good strong scaling results for pivot selection, and better scaling results for core construction. The weak scaling results follow the same behavior seen in the $d = 6$ Hilbert tensor, as well as the 2d2v Maxwellian tensor. Furthermore,~\cref{fig:Max_6D_Strong_Weak_Err_Acc} (Bottom left) includes the error plot, which indicates that our algorithm can obtain a close approximation, even with high dimension and a real-world example. Also included in~\cref{fig:Max_6D_Strong_Weak_Err_Acc} (Bottom right), we have a plot of percentage of true tensor data required for construction. The results are similar to those in the Hilbert case and shows that this method maintains a low storage requirement for TT core construction. In testing, this local storage requirement ranges from $1.47$ Gb for 1 MPI rank down to $4.51\times 10^{-2}$ Gb for 64 MPI ranks to run the full algorithm.

\subsection{Comparisons with dimension parallel and torchTT}
\label{sec:Comparison}
\reva{
In the first part of this subsection, we make a comparison of timings for pivot selection of our proposed algorithm with our own implementation of the dimension parallel TT-Cross algorithm in \cite{dolgov2020parallel}. The main idea of the dimension parallel algorithm is to select pivots across all dimensions simultaneously, and then communicate indices to neighboring dimensions to update superblocks for the next iteration of searching. The timings reported are in the fully parallel regime of both algorithms for all tensors tested in previous sections. Note that the compression of both algorithms is the same as each is prescribed the same core ranks. Furthermore, the dimension parallel algorithm reported here is only in parallel for index selection, and not core construction. 

As seen in~\cref{tab:subtensorvsdimension}, in both Hilbert tensors tested in~\cref{sec:HilbertTensor} as well as the Maxwellian tensors in~\cref{sec:Maxwellian}, the subtensor parallel algorithm using 64 MPI ranks has significantly smaller run time for index selection compared to the dimension parallel algorithm. These results are expected as the mode size of the tensors are quite large, and the subtensor parallel algorithm uses more computing resources.

\begin{table}[h!]
    \centering
    \begin{tabular}{|c|c|c|c|c|}
    \hline
        Alg./Tensor & Hilbert $3d$ & Hilbert $6d$ & Maxwellian $4d$ & Maxwellian $6d$  \\
        \hline
        Subtensor & $99.99$& $260.72$& $147.34$& $184.16$ \\
        Dimension & $2482.75$& $950.43$& $1743.47$ & $538.54$ \\
        \hline
    \end{tabular}
    \caption{Timings in seconds for pivot selection of subtensor parallel compared to dimension parallel in a fully parallel regime. The mode sizes for the various tests are: Hilbert 3$d$: (2000,2000,2000), Hilbert 6$d$: (300,300,300,300,300,300), Maxwellian 4$d$: (2000,1000,2000,1000), Maxwellian 6$d$: (800,400,800,400,800,400).}
    \label{tab:subtensorvsdimension}
\end{table}

We also compare timings of subtensor parallel and dimension parallel on a Hilbert tensor as dimension increases in two more settings: (1) the number of entries is held approximately constant, and (2) the mode size is fixed. For the first case, the corresponding mode sizes $n$ to dimension $d$ are: $n = 1400$ for $d = 5$, $n = 419$ for $d = 6$, $n = 177$ for $d = 7$, $n = 93$ for $d = 8$, $n = 56$ for $d = 9$, and $n = 37$ for $d = 10$, so that the all of the test cases have roughly $5.3e15$ entries. The choice of total size $5.3e13$ is selected as to run the largest test without running into memory request errors. In the second case, we fix $n = 100$ for all cases of dimension $d$.

\begin{figure}
    \centering
    \includegraphics[width=0.5\linewidth]{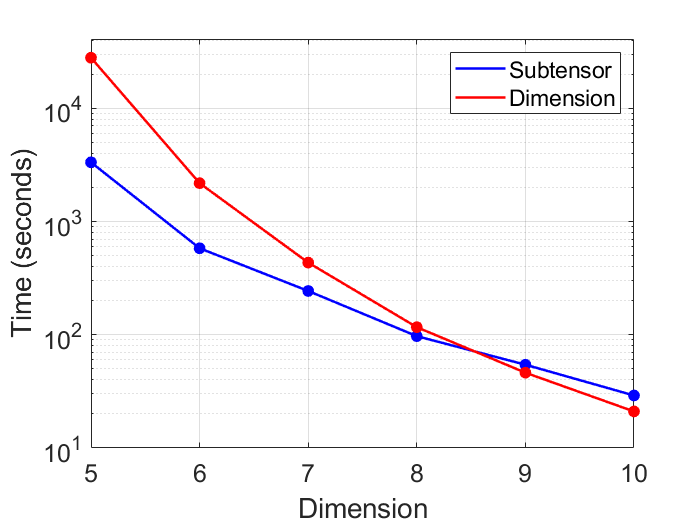}\includegraphics[width=0.5\textwidth]{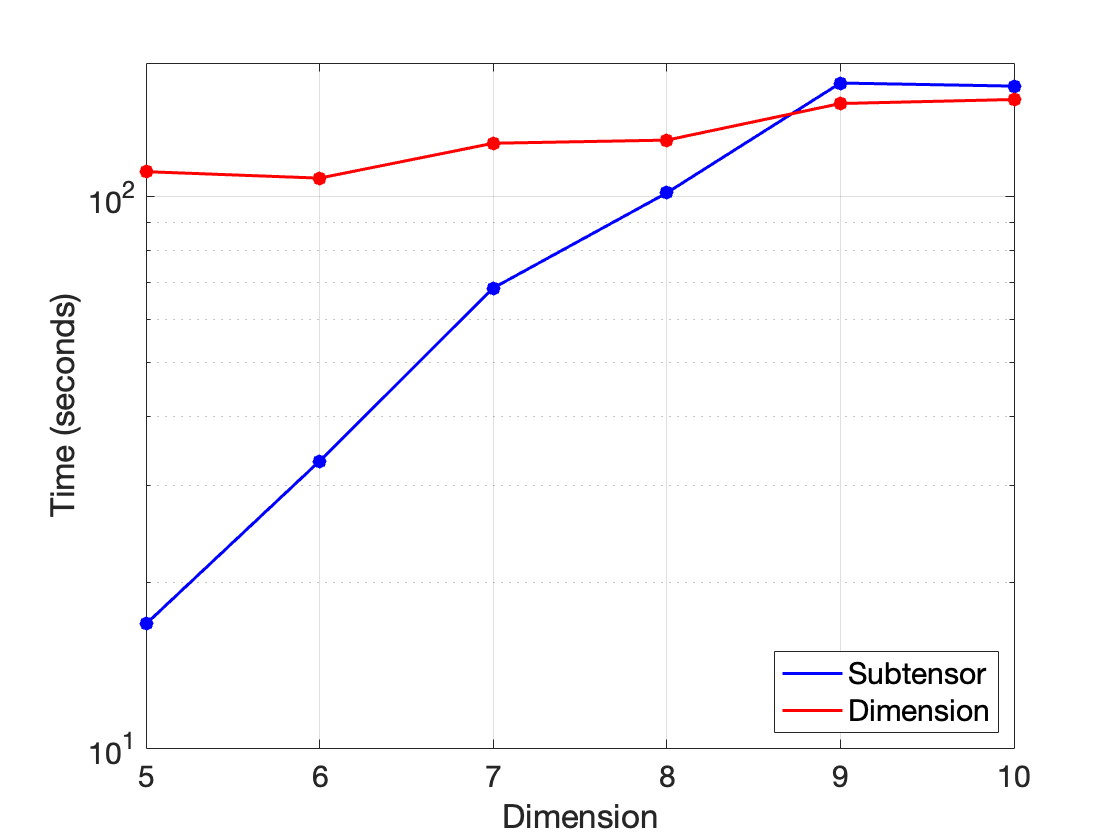}
    \caption{Timings in seconds for index selection using subtensor parallel and dimension parallel. Left: Constant tensor size with varying mode size. Right: Constant mode size with varying dimension.}
    \label{fig:compareplots}
\end{figure}

As shown in both plots of~\cref{fig:compareplots}, for the case of lower dimensions we observe a significant reduction in run time for index selection. Both tests also demonstrate that when dimension grows, the subtensor parallel algorithm approaches the run time of the dimension parallel algorithm. This behavior is expected as the dimension parallel algorithm is designed for tensors with large dimensionality and relatively smaller mode size. Therefore, as the fundamental framework of the subtensor and dimension parallel algorithms are distinct, we are unsurprised to see a transition point where the two perform similarly. 

In the last part of this section, we compare the subtensor parallel algorithm with the publicly available python package torchTT~\cite{torchtt2024}, which contains a TT-Cross function. In their implementation, they do not use a greedy search such as~\cref{alg:MatCross}, but rather a maxvol search \cite{savostyanov2014quasioptimality}. For our implementation of torchTT, we use a pure python execution that runs in serial. In this comparison, we only test on a $3d$ Hilbert tensor, as the $6d$ Hilbert tensor encounters insufficient memory errors when running torchTT. For $4d$ and $6d$ Maxwellian tensors, the observed time to run torchTT on these tensors significantly exceeds the times achieved by the subtensor algorithms. In a small test case with mode size $n = 20$, torchTT takes roughly $1200$ seconds to compress a $6d$ Maxwellian. For this reason, we choose to not attempt a comparison for the other tensors considered in this section. 

\begin{table}[h!]
    \centering
    \begin{tabular}{|c|c|c|c|c|c|c|}
    \hline
        \multirow{2}*{Alg/$n$} & \multicolumn{2}{|c|}{Subtensor} & \multicolumn{2}{|c|}{TT-ACA} & \multicolumn{2}{|c|}{torchTT}  \\
        \cline{2-7} 
        & Time (s) & Error & Time (s) & Error & Time (s) & Error  \\
        \hline
        $250$ & 1.71& 1.27e-06 & 39.02& 1.40e-06& 215.34&6.35e-07 \\
        \hline
        $500$ & 5.06 & 9.99e-07 & 168.09& 5.45e-07& 860.31&2.93e-06 \\
        \hline
        $1000$ &16.07 & 9.07e-07 & 799.42& 1.33e-06& 3560.26&2.27e-06 \\
        \hline
        $2000$ & 60.54 & 8.64e-07 & 3605.90& 3.76e-07& 25159.67&3.19e-06 \\
        \hline
    \end{tabular}
    \caption{Timings in seconds and error for full approximation using subtensor parallel with 64 MPI ranks, TT-ACA (subtensor parallel with one process), and torchTT for a $3d$ Hilbert tensor with mode size $n$.}
    \label{tab:subtensorvstorchttHilbert}
\end{table}

As shown in~\cref{tab:subtensorvstorchttHilbert}, the time taken by the subtensor algorithm is a substantial reduction compared to the times taken by torchTT. It is suspected that one contributing factor to the increase in time is the selection procedure. While maxvol is a robust algorithm for index selection, it is more computationally intensive compared to~\cref{alg:MatCross}, and the accuracy comparison in~\cref{tab:subtensorvstorchttHilbert} does not show much difference between the two schemes. If we assume that torchTT is implemented with perfect parallel scaling, then with the same number of processes used for the subtensor results of~\cref{tab:subtensorvstorchttHilbert}, torchTT would still report times ranging from $3.36$ seconds to $393.12$ seconds. Therefore, even in a perfect scenario, we still observe substantial reduction in computational time.}

\section{Conclusions}
\label{sec:Conclusions}

We introduce a new distributed memory subtensor parallel algorithm for constructing a TT cross approximation for a given tensor. This construction allows for local computations with little communication requirements to obtain TT cores at a global level, while maintaining the interpolation property of a standard TT cross algorithm. This algorithm relies on the efficient update method for pivot selection using the material of~\cref{sec:itersubMat1}, which ensures the selection of the optimal global pivot for each dimension. Furthermore, we utilize multiple process grids of~\cref{sec:TTcrossAlg} combined with an alternate recursive update formula~\cref{sec:TTcrossCore} to construct all TT core simultaneously, with each core constructed in a distributed memory framework. Our presented numerical results demonstrate results ranging from good to optimal scaling for both 3d and 6d Hilbert tensors as well as real world 4d and 6d Maxwellian datasets. \reva{Furthermore, in comparisons with a dimension parallel algorithm for index selection, as well as publicly available package torchTT, we observe significantly lower run time in all test cases.}

There are several directions we will work on in the future. First, in the scaling results in~\cref{sec:NumericalExamples}, we see curves behave beyond optimal a couple of times, and our immediate future work is to analyze and understand these behaviors, which can assist us in building faster and more stable implementations. Furthermore, as subtensors are straightforward vessels for distributed memory parallel tensor algorithms, we want to develop a subtensor parallel TT decomposition framework for all suitable TT factorization algorithms, with thorough analysis of computational complexity, storage costs, and communication studies. Currently, we know this framework can include the adaptive TT cross described in this article, and randomized SVD based TT sketching. There are a few potential algorithm candidates on our radar for this framework, and we will start with column-pivoted QR (CPQR) and LU based TT cross.

\printbibliography

\end{document}